\documentstyle[amscd,amssymb,verbatim,diagrams,12pt]{amsart}
\newcommand{\pa}{\partial}

\newcommand{\ev}{\operatorname{ev}}

\renewcommand{\mod}{\operatorname{mod}}

\newcommand{\Cone}{\operatorname{Cone}}

\newcommand{\OO}{{\cal O}}

\newcommand{\G}{{\Bbb G}}
\newcommand{\A}{{\Bbb A}}

\newcommand{\mg}{{\frak m}}
\newcommand{\hra}{\hookrightarrow}
\newcommand{\lan}{\langle}
\newcommand{\ran}{\rangle}

\newcommand{\CC}{{\cal C}}

\newcommand{\supp}{\operatorname{supp}}
\newcommand{\Spec}{\operatorname{Spec}}

\renewcommand{\P}{{\Bbb P}}

\newcommand{\ga}{\gamma}
\newcommand{\de}{\delta}
\newcommand{\eps}{\epsilon}

\renewcommand{\ker}{\operatorname{ker}}
\newcommand{\im}{\operatorname{im}}

\numberwithin{equation}{subsection}

\newtheorem{thm}{Theorem}[subsection]
\newtheorem{prop}[thm]{Proposition}
\newtheorem{lem}[thm]{Lemma}
\newtheorem{cor}[thm]{Corollary}
{  \theoremstyle{definition}

\newtheorem{ex}[thm]{Example}

\newtheorem{rem}[thm]{Remark}
\newtheorem{rems}[thm]{Remarks}
}

\newcommand{\Pf}{\noindent {\it Proof}}
\newcommand{\id}{\operatorname{id}}

\newcommand{\ov}{\overline}

\newcommand{\rk}{\operatorname{rk}}

\renewcommand{\AA}{{\cal A}}

\newcommand{\MM}{{\cal M}}
\newcommand{\TT}{{\cal T}}
\newcommand{\UU}{{\cal U}}

\newcommand{\LL}{{\cal L}}

\newcommand{\Om}{\Omega}

\newcommand{\Hom}{\operatorname{Hom}}

\newcommand{\Ext}{\operatorname{Ext}}

\newcommand{\Res}{\operatorname{Res}}

\renewcommand{\a}{\alpha}
\renewcommand{\b}{\beta}
\newcommand{\om}{\omega}

\newcommand{\la}{\lambda}
\newcommand{\th}{\theta}
\newcommand{\C}{{\Bbb C}}

\newcommand{\Z}{{\Bbb Z}}

\newcommand{\La}{\Lambda}
\newcommand{\Ga}{\Gamma}

\newcommand{\bM}{{\bf M}}
\newcommand{\wt}{\widetilde}
\newcommand{\ot}{\otimes}

\newcommand{\sub}{\subset}
\newcommand{\ed}{\qed\vspace{3mm}}

\newcommand{\cha}{\operatorname{char}}

\newsymbol\k 207C

\newcommand{\conv}{\operatorname{conv}}

\title{$A_{\infty}$-algebras associated with curves and rational functions on $\MM_{g,g}$. I.}
\author{Robert Fisette and Alexander Polishchuk}
\subjclass[2010]{Primary 14F05; Secondary 16E45, 14H10, 55S30}
\keywords{$A_\infty$-structure, Massey products, derived category, Hochschild cohomology, algebraic curve}
\thanks{Supported in part by the NSF grant DMS-1001364}

\begin{document}

\begin{abstract} We consider the natural $A_{\infty}$-structure on the $\Ext$-algebra $\Ext^*(G,G)$
associated with the coherent sheaf $G=\OO_C\oplus\OO_{p_1}\oplus\ldots\oplus\OO_{p_n}$
on a smooth projective curve $C$, where $p_1,\ldots,p_n\in C$ are distinct points. 
We study the homotopy class of the product $m_3$. 
Assuming that $h^0(p_1+\ldots+p_n)=1$ we prove that $m_3$ is homotopic to zero if and only if
$C$ is hyperelliptic and the points $p_i$ are Weierstrass points. In the latter case we show
that $m_4$ is not homotopic to zero, provided the genus of $C$ is $>1$.
In the case $n=g$ we prove
that the $A_{\infty}$-structure is determined uniquely (up to homotopy) by the products
$m_i$ with $i\le 6$. Also, in this case we study the rational map $\MM_{g,g}\to\A^{g^2-2g}$
associated with the homotopy class of $m_3$. We prove that for $g\ge 6$ it is birational onto
its image, while for $g\le 5$ it is dominant. We also give an interpretation of this
map in terms of tangents to $C$ in the canonical embedding and in the projective embedding given by the linear series $|2(p_1+\ldots+p_g)|$.
\end{abstract}

\maketitle


\centerline{\sc Introduction}

Let $C$ be a smooth projective curve of genus $g$ over an algebraically closed field $\k$.
With any generator $G$ of the derived category $D^b(C)$ of coherent sheaves on $C$
one can associate an $A_{\infty}$-algebra of endomorphisms of $G$, which is basically
the Ext-algebra $\Ext^*(G,G)$ equipped with higher operations defined uniquely
up to homotopy. More precisely, this construction uses a dg-enhancement
of $D^b(C)$ and applies to it the homological perturbation theory 
developed originally in \cite{GS}, \cite{GLS}, \cite{Kadeishvili}, with explicit formulas
given in \cite{Merkulov}, \cite{KS}.
Furthermore, this $A_{\infty}$-algebra determines the derived category $D^b(C)$ (see
\cite[Thm.\ 3.1]{Keller}), and hence the
curve $C$ (at least, if either $\cha(\k)=0$ or $g\neq 1$; see \cite{HVdB}).

One of the possible choices of a generator of $D^b(C)$ is $G=\OO_C\oplus L$,
where $L=\OO_C(p)$ for some point $p\in C$ (the subcategory generated by $G$ contains
$\OO_C(np)$ for any $n>0$, so $G$ is a generator by \cite[Thm.\ 4]{Orlov}). 
In the case of an elliptic curve 
the corresponding $A_{\infty}$-algebra
was explicitly computed in \cite{P-ell} (assuming $\k=\C$).
Note also that in this case there exists an autoequivalence of $D^b(C)$ sending
$G$ to $\OO_C\oplus\OO_p$. Also, in the genus $1$ case Lekili and Perutz studied in \cite{LP} 
homotopy classes of minimal $A_{\infty}$-structures on $\Ext^*(G,G)$
extending the natural double product.
Their results imply that all nontrivial homotopy classes of such 
$A_{\infty}$-structures arise either from
elliptic curves or from the nodal plane cubic (see \cite[Prop.\ 9]{LP}).
Also, any such $A_{\infty}$-structure is finitely determined, i.e., determined up to homotopy
by a finite number of the products $m_i$ (actually by $m_i$ with $i\le 8$).

In this paper we consider a partial extension of this picture to higher genus curves and to the case
of generators of $D^b(C)$ of the form
\begin{equation}\label{G-generator-eq}
G=\OO_C\oplus\OO_{p_1}\oplus\ldots\oplus\OO_{p_n},
\end{equation}
where $p_1,\ldots,p_n$ are $n$ distinct points on $C$ ($n\ge 1$),
such that $h^0(p_1+\ldots+p_n)=1$ (in particular, $n\le g$). 
We would like to study the $A_{\infty}$-structure
on the corresponding $\Ext$-algebra
$$E=E_{g,n}=\Ext^*(G,G),$$
which, depends only on $n$ and $g$ as an associative algebra (however, the higher products
depend on $(C,p_1,\ldots,p_n)$).
In the case $n=g$ this $\Ext$-algebra looks particularly nice:
it is generated over the $(g+1)$-dimensional subalgebra spanned by the natural
idempotents in $\Hom(G,G)$, by the one-dimensional spaces
$\Hom(\OO_C,\OO_{p_i})$ and $\Ext^1(\OO_{p_i},\OO_C)$. Furthermore, the defining relations
between these generators are monomial (see \eqref{relations-eq}).

In general, it follows from the work of To\"en (see \cite[Cor.\ 1.3]{Toen}) that for any generator $G$
of the derived category of coherent sheaves on a smooth projective variety the natural $A_\infty$-structure on
$\Ext^*(G,G)$ is determined up to homotopy by a finite number of the products $m_i$ (in the class
of $A_\infty$-structures that are smooth and proper). It turns out that in our case a stronger statement holds:
any minimal $A_{\infty}$-structure on $E_{g,g}$ is finitely determined. More precisely,
it was proved in \cite{Fisette-thesis} that any minimal $A_\infty$-structure on $E_{g,g}$ is
determined up to homotopy by $m_i$ with $i\le 6$. This follows from the vanishing of certain
graded components of the Hochschild cohomology of $E_{g,g}$ 
(see Theorem \ref{finiteness-thm}).
We give a simpler proof of this vanishing using a minimal resolution of
$E_{g,g}$ from \cite{Bardzell}.
On the other hand, we show that the same 
vanishing does not hold for the algebra $E_{g,n}$ if $n<g$,
and in this case there exist $A_{\infty}$-structures on $E_{g,n}$ that are not determined
by any fixed finite number of $m_i$ (with a possible exception of the case $g=2$, $n=1$;
see Remark \ref{finiteness-rem}.2). 

We also consider the following basic question about the $A_{\infty}$-structure on 
$E_{g,n}$ coming from $(C,p_1,\ldots,p_n)$: whether it is equivalent to the one with $m_3=0$.
We prove that this holds if and only if $C$ is hyperelliptic and the points $p_1,\ldots,p_n$
are Weierstrass points (see Theorem \ref{m3-m4-thm}). Furthermore, we also show that
if $g>1$ then either $m_3$ or $m_4$ is always nontrivial.
The main point in the proof is that the Hochschild cohomology class given by the triple product $m_3$
can be recovered from the triple Massey products for the complexes
\begin{equation}\label{main-triple-Mas-prod-eq}
\OO\rTo{}\OO_{p_i}\rTo{[1]}\OO_{p_i}\rTo{[1]}\OO
\end{equation}
(see Proposition \ref{Hoch-prop} and Section \ref{3Mas-sec}). 
In the hyperelliptic case we also study a certain quadruple
Massey product and use \cite[Thm.\ 3.1]{LPWZ} to connect it with $m_4$.

In the case $n=g$ we compute the triple Massey products \eqref{main-triple-Mas-prod-eq} in terms of canonical
rational sections of some natural line bundles on the moduli spaces
$\MM_{g,g}$ of curves with $g$ marked points. Considering rational monomials of these sections,
we get $g^2-2g$ rational functions on $\MM_{g,g}$, i.e., a rational
map
\begin{equation}\label{moduli-rat-map}
\ov{\a}:\MM_{g,g}\to\A^{g^2-2g}
\end{equation}
(see Section \ref{rat-fun-sec}).
Assuming that the characteristic is zero, 
we prove that for $g\ge 6$ this map is birational onto its image (see Theorem
\ref{birational-thm}), while for $g\le 5$ it is dominant (see Theorem \ref{dominant-thm}).
The main idea in the proof of the former result is to reconstruct a curve $C$ from the
multiplication table between certain rational functions with prescribed polar parts
at $p_1,\ldots,p_g\in C$ (see Section \ref{recon-sec}). We also observe
that the above rational map extends to stable curves and make explicit calculations for
rational irreducible nodal curves (see Section \ref{degen-sec})
To prove dominance for $g\le 5$
we first calculate the tangent map (see Section \ref{tangent-map-sec}). Then we again use explicit calculations for rational
nodal curves.

It is interesting to note that 
our triple Massey products \eqref{main-triple-Mas-prod-eq}
have a nice geometric interpretation: they record positions of the tangent lines to $C$ at $p_i$
in the canonical embedding, as well as, for $n=g$, of the tangent lines to $C$ at $p_i$ in
the projective embedding given by the linear system $|2(p_1+\ldots+p_g)|$. 
Equivalently, they can be related to the Wahl maps (defined in \cite{Wahl}) for the line bundles
$\om_C$ and $\OO(2(p_1+\ldots+p_g))$ evaluated at the marked points 
(see Section \ref{geom-int-sec}). 

The interest in characterizing $A_{\infty}$-algebras of the form $\Ext^*(G,G)$ is motivated
by the homological mirror symmetry conjecture, extended to non-Calabi-Yau manifolds
(see \cite{KKOY}). Note that one knows the homological mirror correspondence
involving a higher genus curve on the symplectic side and a Landau-Ginzburg model on the B-side
due to the work of Seidel and Efimov \cite{Seidel}, \cite{Efimov}. 
The other half of the correspondence for the same
mirror pair should involve the derived category $D^b(C)$, governed by the 
$A_{\infty}$-algebra $\Ext^*(G,G)$. Thus, finding a
characterization of $A_{\infty}$-structures on $E_{g,g}$ arising from curves would be a step
towards establishing such a correspondence. In a sequel to this paper we will compute explicitly
the higher products on $E_{g,g}$ arising from curves and will study the normal forms of
arbitrary $A_{\infty}$-structures on $E_{g,g}$.

The paper is organized as follows. In Section \ref{Hoch-sec} we perform the calculation
of the relevant Hochschild cohomology of the algebras $E_{g,n}$ (mostly for $n=g$).
Section \ref{Massey-sec} is devoted to Massey products. Here we compute the triple
Massey products governing the Hochschild cohomology class of $m_3$ on $E_{g,n}$,
and a certain quadruple Massey product related to $m_4$. This allows us
to characterize geometrically the vanishing of $m_3$ (see Theorem \ref{m3-m4-thm}).
In Section \ref{rat-sec} we study the Massey products \eqref{main-triple-Mas-prod-eq} globally over the moduli space of curves and show how they lead to the rational map \eqref{moduli-rat-map}.
Also, in Section \ref{geom-int-sec} we discuss the connection with the tangent lines to $C$ in the canonical embedding and in the embedding given by $|2(p_1+\ldots+p_g)|$ and with the 
corresponding Wahl maps.
In Section \ref{recon-sec} we prove that \eqref{moduli-rat-map} is birational onto its image for
$g\ge 6$. Finally, in Section \ref{tangent-map-sec} we compute the tangent map
to \eqref{moduli-rat-map} and show that it is dominant for $g\le 5$.
The Appendix contains GAP codes that we used to carry out explicit calculations 
needed for some proofs.

\medskip

\noindent
{\it Notation and conventions}. We work over a fixed ground field $\k$, which is assumed to be algebraically
closed whenever we discuss geometry. By a curve we mean a projective connected curve over $\k$.
By a divisor on a (not necessarily smooth) curve $C$ we always mean
a divisor supported on the smooth part of $C$. For such a divisor $D$ we use
the notation $h^i(D)=\dim_\k H^i(C,\OO(D))$ for $i=0,1$.
We use the similar notation $h^i(L)$ for a line bundle $L$.
In a triangulated category we denote
$\Hom^n(X,Y):=\Hom(X,Y[n])$ for $n\in\Z$. We also depict elements of $\Hom^n(X,Y)$ by
arrows $X\rTo{[n]} Y$. For a morphism $f:X\to Y$ we often denote the morphism $f[n]:X[n]\to Y[n]$ 
simply by $f$.
For a dg-category $\CC$ we denote the differentials on the $\Hom$-spaces by $\pa$.
We denote by $H^*(\CC)$ (resp., $H^0(\CC)$) the category obtained by passing to cohomology
(resp., $0$th cohomology) in $\Hom$-spaces of $\CC$.
For dg $\CC$-modules $M,N$ we denote by $\Hom_\CC(M,N)$ the
space of morphisms in the dg-category of dg $\CC$-modules.
All our $A_\infty$-structures are assumed to be strictly unital.
For a vector space $V$ with a basis $B$, an element $b\in B$, and an element $w\in W$ in another vector space, we denote by $[b]^*w$ the linear map $V\to W$ sending $b$ to $w$ and 
$B\setminus b$ to zero.
For a line bundle or a $1$-dimensional vector space $L$ we often abbreviate
$L^{\ot n}$ as $L^n$.

\medskip

\noindent
{\it Acknowledgments}. The second author would like to thank Michael Kapranov and Ravi Vakil
for useful discussions.

\section{Hochschild cohomology}\label{Hoch-sec}

We refer to \cite{Keller-intro} for an introduction to $A_{\infty}$-algebras.
Recall that an $A_{\infty}$-structure $(m_i)$ on a vector space $A$ is called {\it minimal} if $m_1=0$.
In this case $m_2$ equips $A$ with a structure of a graded associative algebra. Thus, fixing $m_2$ we can
talk about minimal $A_{\infty}$-structures on a graded algebra $A$.
It is well known that equivalence classes of such $A_{\infty}$-structures on $A$
are controlled by the Hochschild cohomology
$HH^*(A)=H^*(A,A)$. In particular, if we have two such structures $(m_i)$ and $(m'_i)$ with $m_i=m'_i$
for $i<n$ then $m'_n-m_n$ is a Hochschild $n$-cocycle of internal degree $2-n$, 
whose triviality means that the structure $(m'_i)$
can be changed by a homotopy in such a way that $m_i=m'_i$ for $i\le n$ (see \cite[Lem.\ 2.2]{P-ext}). 
Let us denote by $HH^i(A)_j$ the component of the $i$th Hochschild cohomology group
of internal degree $j$.
We deduce that the vanishing of the Hochschild cohomology $HH^i(A)_{2-i}$ for all $i>N$ implies that
any minimal $A_{\infty}$-structure on $A$ is determined by $(m_i)$ with $i\le N$ up to homotopy.

We would like to apply these principles to the $\Ext$-algebra $E=E_{g,n}$ (described explicitly below).
In the case $n=g$ the relevant Hochschild cohomology was studied in \cite{Fisette-thesis}.
Here we present two results of this study: an explicit description of $HH^3(E)_{-1}$
and the vanishing of $HH^i(E,E)_{2-i}$ for large $i$ (see Proposition \ref{Hoch-prop}
and Theorem \ref{finiteness-thm}(i) below). In addition, we will show that
the latter property does not hold if $n<g$.

\subsection{Algebras $E_{g,n}$}\label{Egn-sec}

Let $C$ be a projective curve over $\k$ of arithmetic genus $g$, and let $p_1,\ldots,p_n$ be distinct smooth points
on $C$ such that $h^0(p_1+\ldots+p_n)=1$. Then from the short exact sequence
$$0\to \OO_C\to\OO_C(p_1+\ldots+p_n)\to\OO_C(p_1+\ldots+p_n)/\OO_C\to 0$$
we get the boundary homomorphism
$$\bigoplus_{i=1}^n H^0(C,\OO(p_i)/\OO)\simeq H^0(C,\OO(p_1+\ldots+p_n)/\OO)\to 
H^1(C,\OO),$$
which is an embedding, since the map $H^0(C,\OO)\to H^0(C,\OO(p_1+\ldots+p_n))$ is
an isomorphism by our assumption. 

Let $A_i\in\Hom(\OO_C,\OO_{p_i})$ and 
$B_i\in\Ext^1(\OO_{p_i},\OO_C)$ be generators of these one-dimensional spaces.
Then 
$$Y_i=A_i\circ B_i\in\Ext^1(\OO_{p_i},\OO_{p_i})$$
is a generator of $\Ext^1(\OO_{p_i},\OO_{p_i})$, and
the elements
$$X_i=B_i\circ A_i\in\Ext^1(\OO_C,\OO_C)=H^1(\OO_C)$$
for $i=1,\ldots,n$ are linearly independent.
In the case $n<g$ we extend these to a basis $(X_1,\ldots,X_g)$ of $H^1(\OO_C)$.

Thus, the algebra $E_{g,n}=\Ext^*(G,G)$, where $G$ is given by \eqref{G-generator-eq}, 
has the $\k$-basis 
\begin{align*}
&e_{\OO}:=\id_{\OO}, \ e_{\OO_{p_i}}:=\id_{\OO_{p_i}}, \ A_i, \ B_i, \ Y_i, \ i=1,\ldots,n;\\
&X_j, \ j=1,\ldots,g.
\end{align*}
The only nontrivial products in $E_{g,n}$ are the obvious relations involving the idempotents
$e_\OO$ and $e_{\OO_{p_i}}$, as well as the relations $A_iB_i=Y_i$ and
$B_i A_i=X_i$ for $i=1,\ldots,n$. In particular, this algebra does not depend on a specific
curve and points on it.

Note that the algebra $E_{g,n}$ is the quotient algebra of the path algebra of the quiver
$\Ga_{g,n}$ with $n+1$ vertices, marked with $\OO$ and $\OO_{p_i}$, $i=1\ldots,n$.
The arrows in $\Ga_{g,n}$ go in the direction opposite to the direction of morphisms in $D^b(C)$.
Namely, for each $p_i$ we have one arrow of degree $1$ from $\OO$ to $\OO_{p_i}$ and
one arrow of degree $0$ in the opposite direction. In addition, we have $g-n$ loops
of degree $1$ at $\OO$ (that correspond to the generators $X_{n+1},\ldots,X_g$).

We denote by $E_{g,n}^+$ the ideal in $E_{g,n}$ obtained from paths of length $\ge 1$.
In other words, this is the $\k$-subspace spanned by all $A_i$, $B_i$, $Y_i$ and $X_j$.

\subsection{Minimal resolution of $E_{g,g}$}\label{resolution-sec}

Our method of calculating the Hochschild cohomology of $E=E_{g,g}$ is similar to that
used in \cite{LP} for $g=1$. Namely, 
we view $E$ as the quotient of the path algebra $\k[\Ga_{g,g}]$ by the monomial
relations 
\begin{equation}\label{relations-eq}
A_iB_i A_i=B_i A_iB_i=A_iB_j=0
\end{equation}
for $1\le i,j\le g$, $i\neq j$.
Hence, we can use a minimal projective resolution 
$$\ldots\to P_1\to P_0\to E$$ 
over the enveloping algebra $E^e=E\ot E^{op}$ constructed in \cite{Bardzell}. 
Let us recall this construction.

For every pair of vertices $v,v'$ in the quiver we have the projective $E-E$-bimodule
$$P_{v,v'}:=Ee_v\ot e_{v'}E,$$ 
where $e_v$ is the idempotent in $E$
corresponding to $v$. For a path $p$ in $\Ga_{g,g}$ let $v$ and $v'$ be vertices such that
$e_vpe_{v'}=p$ in $\k[\Ga_{g,g}]$. Then we call $P_{v,v'}$
the {\it projective bimodule generated by} $[p]$, and denote its elements by $x[p]y$, where 
$x\in Ee_v$, $y\in e_{v'}E$. We define the projective bimodule generated by a collection of paths
as the direct sum of the projective bimodules generated by each path.

The $E-E$-bimodules in our minimal resolution are defined as follows: 
$P_0=E^e$, and for $j>0$ we define
$P_j$ as the projective bimodule generated by the set $AP(j)$ of paths in $\Ga_{g,g}$,
defined by the following recursive procedure\footnote{We specialize a more
general procedure from \cite{Bardzell} to our situation}.
By definition, $AP(1)$ consists of all paths of length $1$, i.e., of $A_i$ and $B_i$
($i=1,\ldots,g$), while $AP(2)$ is exactly the set $R$ of generating relations, namely of the paths
in \eqref{relations-eq}. Next, $AP(3)$ consists of paths ``linking"
pairs from $R$, namely, 
$$AP(3)=\{(A_iB_i)^2, (B_i A_i)^2, A_iB_i A_iB_j, A_iB_j A_jB_j
\ |\ 1\le i,j\le g, i\neq j\}.$$
Let us denote by $S$ the set of nonempty proper subwords in $R$. Thus,
$$S=\{A_iB_i, B_i A_i, A_i, B_i \ |\ 1\le i\le g\}.
$$
Note that each path in $AP(3)$ has the form $p=sr$, where $r\in R$ and $s\in S$. 
Similarly, every path in $AP(j)$ will be of the form $p=sp'$, where $p'\in AP(j-1)$ and
$s\in S$. By definition, for $j\ge 3$, $AP(j+1)$ is obtained by taking all paths $p=sp'\in AP(j)$
(where $p'\in AP(j-1), s\in S$)
and replacing $s$ either by an element of $R$ ending with $s$ or, in the case $s=B_i A_i$,
by $A_jB_i A_i$. For example,
\begin{equation}\label{AP4-eq}
AP(4)=\{(A_iB_i)^3, (B_i A_i)^3, (B_i A_i)^2B_j, A_iB_i A_iB_j A_jB_j,
A_i(B_j A_j)^2 \ |\ 1\le i,j\le g, i\neq j\}.
\end{equation}
For a path $p=sp'\in AP(j)$ with $p'\in AP(j-1)$ and $s\in S$, we call $s$ the {\it head of $p$}.

For an arrow $a$ in the quiver $\Ga_{g,g}$ we denote by $s(a)$ and $t(a)$ the source and
the target of $a$ (these are vertices of $\Ga_{g,g}$). 
The first two of the differentials $d_j:P_j\to P_{j-1}$ are described as follows:
$$d_1: [a]\mapsto e_{s(a)}\ot a-a\ot e_{t(a)},$$
$$d_2: [a_1a_2\ldots] \mapsto [a_1]a_2\ldots+a_1[a_2]\ldots+\ldots,$$
where $a\in AP(1)$, $a_1a_2\ldots\in AP(2)$.
For odd $j>2$ the differential is
$$d_j: [p]\mapsto s'[p']-[p'']s'',$$
where we write $p\in AP(j)$ in the form $p=s'p'=p''s''$ with $p',p''\in AP(j-1)$ and $s,s'\in S$.
For even $j>2$ the differential is
$$d_j: [p]\mapsto \sum s_1[p']s_2,$$
where $p\in AP(j)$ and the sum is over all decompositions $p=s_1p's_2$ with $p'\in AP(j-1)$.
It is shown in \cite[Thm.\ 4.1]{Bardzell} that we get in this way a minimal projective resolution
of $E$ over $E^e$.

\begin{lem}\label{Pj-deg-lem} 
(a) The maximal internal degree of the generators of $P_j$ is equal
to 
$$h(j):=\begin{cases} j-[j/4]-1, &j\equiv -1 \mod (4);\\ j-[j/4], &\text{otherwise}.\end{cases}
$$

\noindent
(b) The maximal internal degree of the generators of $P_{10}$ (resp., $P_9$)
that end with $A_i$ is equal to $7$ (resp., $6$).
\end{lem}

\Pf . (a) For each $s\in S$ let us denote by $a_j(s)$ the maximal degree of a word in $AP(j)$
with the head $s$ (where $\deg A_i=0$, $\deg B_i=1$). Then from the definition of $AP(j)$
we get that $a_3(s)=2$ for all $s$, as well as the recursive formulas
\begin{align*}
&a_{j+1}(AB)=a_j(A)+1, \\
&a_{j+1}(B)=a_j(AB)+1,\\
&a_{j+1}(BA)=a_j(B)+1,\\
&a_{j+1}(A)=\max(a_j(B),a_j(BA))
\end{align*}
for $j\ge 3$. Here we omit indices with $A$ and $B$ since the value of $a_j(\cdot)$
does not depend on them.
Now it easy to check by induction that
$$a_j(BA)=h(j), \ a_j(B)=h(j+1)-1, \ a_j(AB)=h(j+2)-2, \ a_j(A)=h(j+3)-3,$$
which implies the assertion since $h(j+1)\le h(j)+1$.

\noindent
(b) For $s\in S$ let us denote by $b_j(s)$ the maximal degree 
of a word in $AP(j)$ that has the head $s$ and ends with $A_i$
(in the case when there are no such words we set $b_j(s)=-\infty$). These numbers
satisfy the same recursive formulas as the numbers $a_j(s)$. From this we get 
$$b_9(A)=b_9(AB)=b_9(B)=6, b_9(BA)=-\infty,$$
$$b_{10}(A)=6, b_{10}(AB)=b_{10}(B)=b_{10}(BA)=7,$$
which implies our claim.
\ed

\subsection{Calculations}

Hochschild cohomology groups $HH^i(E_{g,g})_{2-i}$ were calculated in \cite{Fisette-thesis}. 
Here, using a minimal projective resolution of $E_{g,g}$ over its enveloping algebra,
we give a different proof of the fact that these groups vanish for large $i$.

\begin{thm}\label{finiteness-thm} 
(i) One has $HH^i(E_{g,g})_{2-i}=0$ for $i>8$. If $g>1$ then
$HH^i(E_{g,g})_{2-i}=0$ for $i>6$.

\noindent
(ii) Assume $1\le n<g$. Then $HH^i(E_{g,n})_{2-i}\neq 0$ for all $i\ge 5$.
\end{thm}

\Pf . (i) We can compute the Hochschild cohomology of $E=E_{g,g}$ using the minimal resolution 
$P_\bullet\to E$ from Section
\ref{resolution-sec}. First, we claim that $\Hom_{E^e}(P_i,E(2-i))=0$ for $i>10$.
Indeed, Lemma \ref{Pj-deg-lem}(a) implies that the internal degrees
of generators of $P_i$ are $<i-2$.  Next, for $i=9$ or $10$ we still claim that
$\Hom_{E^e}(P_i,E(2-i))=0$. Indeed, first we observe that in this case
$a_i(s)\ge i-2$ only when $s$ begins with some $B_k$. Thus, the only possibly nontrivial
morphism $P_i\to E(2-i)$ should send a generator $p\in AP(i)$ of degree $i-2$, 
beginning with $B_k$, to $E_0$. But any nonzero homogeneous element of degree $0$ in $E_0$ 
that begins at the vertex $\OO$, is proportional to $e_{\OO}$, so $p$ has to end with $A_j$.
By Lemma \ref{Pj-deg-lem}(b), this contradicts $p$ being of degree $i-2$.

Now assume that $g>1$. We claim that the maps
$$d_9^*:\Hom_{E^e}(P_8,E(-6))\to\Hom_{E^e}(P_9,E(-6)) \ \text{ and}$$
$$d_8^*:\Hom_{E^e}(P_7,E(-5))\to\Hom_{E^e}(P_8,E(-5))$$
are injective and hence $HH^8(E)_{-6}=HH^7(E)_{-5}=0$.
First, let us analyze the spaces $\Hom_{E^e}(P_8,E(-6))$ and $\Hom_{E^e}(P_7,E(-5))$
using methods of Lemma \ref{Pj-deg-lem}. Since $(P_8)_{>6}=0$, the only nonzero
components of $\Hom_{E^e}(P_8,E(-6))$ correspond to generators of degree $6$ in $P_8$
mapping to $E_0$. Among such generators $p\in AP(8)$ beginning with $B_i$ we are only interested
in those that end with some $A_j$ (otherwise, there is no element in $E_0$ to map $p$ to).
Note that $a_8(A)=b_8(A)=b_8(AB)=5$ and $b_8(B)=-\infty$.
In particular, we only need to consider $p$ with the heads $B_i A_i$ or $A_iB_i$. 
In the former case $p$ should end with some $A_j$,
so from the definition of $AP(\cdot)$ we see that $p=(B_i A_i)^3p'$, where
$p'\in AP(4)$ has the head $B_j A_j$ and ends with some $A_k$.
From the list \eqref{AP4-eq} we conclude that $p'=(B_j A_j)^3$.
On the other hand, in the case when $p$ has the head $A_iB_i$ we have
$p=A_iB_i A_iB_j A_jB_j p'$, where $p'\in AP(4)$ has the head $A_jB_j$ and 
ends with $B_i$ (otherwise there is no element in $E_0$ to map $p$ to). 
This gives either $p'=(A_iB_i)^3$ or 
$p'=A_jB_jA_jB_iA_iB_i$.
Thus, $\Hom_{E^e}(P_8,E(-6))$ has the following basis:
\begin{align}\label{basis-maps-eq}
&\a_1(i)=[(B_iA_i)^6]^*e_{\OO}, \nonumber\\
&\a_2(i,j)=[(B_iA_i)^3(B_jA_j)^3]^*e_{\OO}, \  i\neq j, \nonumber\\
&\a_3(i)=[(A_iB_i)^6]^*e_{\OO_{p_i}}, \nonumber\\
&\a_4(i,j)=[A_iB_iA_i(B_jA_j)^3B_iA_iB_i]^*e_{\OO_{p_i}},  \  i\neq j, 
\end{align}
Here we identify $\Hom_{E^e}(P_8,E(-6))$ with the subspace of graded linear maps from 
the vector space with the basis $AP(8)$ to $E(-6)$ and denote by $[p]^*x$ the linear map
that sends $p\in AP(8)$ to $x$ and sends other basis elements to zero.
To show that the images of the basis elements \eqref{basis-maps-eq}
under $d_9^*$ stay linearly independent
it is enough to give some basis elements $\b_1(i),\ldots \b_4(i,j)$ in $\Hom_{E^e}(P_9,E(-6))$,
such that $\b_m(A)$ appears in $d_9^*(\a_m(A))$ but not in $d_9^*(\a_n(*))$ with $n>m$ and
not in $d_9^*(\a_m(A'))$ with $A'\neq A$. For this purpose we take
\begin{align}
&\b_1(i)=[A_j(B_iA_i)^6]^*A_j, \nonumber\\
&\b_2(i,j)=[A_j(B_iA_i)^3(B_jA_j)^3]^*A_j, \  i\neq j, \nonumber\\
&\b_3(i)=[B_i(A_iB_i)^6]^*B_i, \nonumber\\
&\b_4(i,j)=[(B_iA_i)^2(B_jA_j)^3B_iA_iB_i]^*B_i,  \  i\neq j, \nonumber
\end{align}
where in the first line we choose any $j$, different from $i$.

The proof of injectivity of $d_8^*$ is very similar, where we use the 
basis of $\Hom_{E^e}(P_7,E(-5))$ given by
\begin{align}
&\a_1(i)=[(B_iA_i)^5]^*e_{\OO}, \nonumber\\
&\a_2(i,j)=[(B_jA_j)^3(B_iA_i)^2]^*e_{\OO}, \  i\neq j, \nonumber\\
&\a_3(i,j)=[(B_jA_j)^2(B_iA_i)^3]^*e_{\OO}, \  i\neq j, \nonumber\\
&\a_4(i)=[(A_iB_i)^5]^*e_{\OO_{p_i}}, \nonumber\\
&\a_5(i,j)=[A_iB_iA_i(B_jA_j)^2B_iA_iB_i]^*e_{\OO_{p_i}},  \  i\neq j, \nonumber\\
&\a_6(i,j)=[A_iB_iA_i(B_jA_j)^3B_i]^*e_{\OO_{p_i}},  \  i\neq j, \nonumber\\
&\a_7(i,j)=[A_i(B_jA_j)^3B_iA_iB_i]^*e_{\OO_{p_i}},  \  i\neq j \nonumber
\end{align}
and the following basis elements in $\Hom_{E^e}(P_8,E(-5))$:
\begin{align}
&\b_1(i)=[A_j(B_iA_i)^5]^*A_j, \nonumber\\
&\b_2(i,j)=[A_j(B_jA_j)^3(B_iA_i)^2]^*A_j, \  i\neq j, \nonumber\\
&\b_3(i,j)=[A_i(B_jA_j)^2(B_iA_i)^3]^*A_i, \  i\neq j, \nonumber\\
&\b_4(i)=[(A_iB_i)^6]^*Y_i, \nonumber\\
&\b_5(i,j)=[(B_iA_i)^2(B_jA_j)^2B_iA_iB_i]^*B_i,  \  i\neq j, \nonumber\\
&\b_6(i,j)=[(B_iA_i)^2(B_jA_j)^3B_i]^*B_i,  \  i\neq j, \nonumber\\
&\b_7(i,j)=[A_iB_iA_i(B_jA_j)^3B_iA_iB_i]^*Y_i,  \  i\neq j, \nonumber
\end{align}
where in the first line we choose any $j$, different from $i$.

\noindent
(ii) Let us consider the standard complex $(C^\bullet,\de)$ computing the Hochschild cohomology of
the algebra $E=E_{g,n}$. The basis of $E$ as $R$-bimodule gives us a basis 
of each $C^n$ of the form $[w]^*b$, where $b\in E$ is a basis element and
$w$ is a (composable) word of length $n$ in basis elements in $E_+$.
Let us set $X=X_{n+1}$.
We claim that for $i\ge 5$ the element
\begin{equation}\label{xi-cocycle-eq}
c=[X B_1Y_1A_1X^{i-4}]^*X_1+[XX_1B_1A_1X^{i-4}]^*X_1\in C^i
\end{equation}
is a cocycle giving a nontrivial cohomology class. Indeed, 
$$\de\left([X B_1Y_1A_1X^{i-4}]^*X_1\right)=-\de\left([XX_1B_1A_1X^{i-4}]^*X_1\right)
=-[X B_1A_1B_1A_1X^{i-4}]^*X_1,$$
so $\de(c)=0$. On the other hand, note that the product of any two consecutive letters in
the word $w=X B_1Y_1A_1X^{i-4}$ is zero. Hence,
for any basis element $[w']^*b\in C^{i-1}$, such that
$[w]^*X_1$ appears in $\de([w']^*b)$, we have either $w=X w'$ and $X_1=X b$ or
$w=w'X$ and $X_1=bX$. Since $X_1$ is not divisible by $X$ either on the left or right, this
is impossible.
\ed

\begin{rems}\label{finiteness-rem}
1. In the case $g=1$ the spaces $HH^6(E_{1,1})_{-4}$ and $HH^8(E_{1,1})_{-6}$ 
are one-dimensional.
Furthermore, (assuming $\cha(\k)\neq 2,3$) any minimal $A_{\infty}$-structure on $E_{1,1}$ 
extending the natural $m_2$
is equivalent to the one for which $m_3=m_4=m_5=0$, and the Hochschild classes of
$m_6$ and $m_8$ completely determine the $A_{\infty}$-structure up to an equivalence
(see \cite[Thm.\ 5]{LP}). In the case $g>1$ any $A_{\infty}$-structure is determined by
the products $m_i$ with $i\le 6$. However, the situation is more complicated since
$m_3$ is usually nonzero (see Theorem \ref{m3-m4-thm} below). Furthermore, Theorem \ref{birational-thm}
below implies that an $A_{\infty}$-structure on $E_{g,g}$ arising from a generic curve of genus $g$ with
$g$ points, is determined among such $A_{\infty}$-structures by $m_3$ alone. However, it is not clear whether every
generic $A_{\infty}$-structure on $E_{g,g}$ arises geometrically for $g>1$ (this is true for $g=1$).

\noindent
2. Assume that $1\le n<g$ and $g\ge 3$. Then for any $i\ge 5$ there exists a minimal $A_{\infty}$-structure on $E_{g,n}$ with standard $m_2$, such that $m_i$ gives 
a nontrivial Hochschild cohomology
class, and $m_j=0$ for $j\neq 2,i$. Indeed, we can define $m_i$ by a slight modification of 
formula \eqref{xi-cocycle-eq}:
$$m_i=[X B_1Y_1A_1X^{i-4}]^*X_j+[XX_1B_1A_1X^{i-4}]^*X_j\in C^i$$
for any $j\neq 1, n+1$. Then the $A_{\infty}$-axiom is satisfied for $(m_2,m_i)$. Hence, in this
case an $A_{\infty}$-structure on $E_{g,n}$ is not determined by any fixed finite number of
$(m_i)$. Note that this does not contradict To\"en's result \cite[Cor.\ 1.3]{Toen}, since we do not
impose the condition for our $A_\infty$-algebras to be smooth.
\end{rems}

\begin{prop}\label{Hoch-prop} Assume $\cha(\k)\neq 2$.
Let us associate with a Hochschild $3$-cocycle $c$ on $E_{g,g}$ of internal degree $-1$
the constants $\a_{ij}(c)$ by
$$c(B_i,Y_i,A_i)=\sum_j \a_{ij}(c) X_j.$$
Then the map 
\begin{equation}\label{alpha-map}
\a: c\mapsto (\a_{ij}(c))_{i\neq j}
\end{equation} 
induces an isomorphism
of $HH^3(E_{g,g})_{-1}$ with the space of $g\times g$-matrices with zeros on the diagonal.
\end{prop}

\Pf . Let us set $E=E_{g,g}$.

\noindent
{\bf Step 1}. First, we check that the map $\a$ is well defined, i.e., that it vanishes on boundaries.
Indeed, for a $2$-cochain $h$ of internal degree $-1$ we have 
$h(B_i,Y_i)=\la\cdot B_i$ and $h(Y_i,A_i)=\la'A_i$ for some constants $\la,\la'$.
Hence,
$$(\de h)(B_i,Y_i,A_i)=-h(B_i,Y_i)\cdot A_i-B_i\cdot h(Y_i,A_i)=
-\la B_iA_i-\la' B_iA_i=-(\la'+\la)X_i,$$
so $\a_{ij}(\de h)=0$ for $i\neq j$.

\noindent
{\bf Step 2}. Next, we claim that the map $\a$ is surjective. Indeed, for $i\neq j$ let us 
consider the Hochschild $3$-cochain
$$f_{ij}=[B_iY_iA_i]^*X_j+[B_iA_iX_i]^*X_j.$$
It is easy to check that $f_{ij}$ is a cocycle and that $\a(f_{ij})$ is the elementary matrix $E_{ij}$.

\noindent
{\bf Step 3}. $\dim HH^3(E)_{-1}=g(g-1)$.
Using the minimal $E^e$-resolution $P_\bullet\to E$ from Section \ref{resolution-sec}
we can identify the space $HH^3(E)_{-1}$ with the middle cohomology in
$$\Hom_{E^e}(P_2,E(-1))\rTo{d_3^*} \Hom_{E^e}(P_3,E(-1))\rTo{d_4^*}\Hom_{E^e}(P_4,E(-1)).$$
First, note that $\Hom_{E^e}(P_4,E(-1))=0$. Indeed, the only generators of degree $\le 2$
in $P_4$ correspond to paths $A_jB_iA_iB_iA_i\in AP(4)$, where $i\neq j$,
and there are no elements of degree $1$ in $e_{\OO_{p_j}} E e_{\OO}$.
The space $\Hom_{E^e}(P_3,E(-1))$ has the basis
$$[B_iA_iB_iA_i]^*X_j, \ [A_iB_iA_iB_i]^*Y_i, \ 1\le i,j\le g.$$
On the other hand, the space $\Hom_{E^e}(P_2,E(-1))$ has the basis
$$[A_iB_iA_i]^*A_i, \ [B_iA_iB_i]^*B_i, \ i=1,\ldots,g.$$
Furthermore, the differential $d_3^*$ is the direct sum of $g$ copies of the same differential as
for the $g=1$ case, which is injective provided $\cha(\k)\neq 2$ (see the proof of \cite[Thm. 4]{LP}).
Hence, the dimension of the cohomology is $(g^2+g)-2g=g^2-g$ as claimed.
\ed

\section{Massey products}\label{Massey-sec}

\subsection{Massey products for dg categories}

Let $(A,\pa)$ be a dg-algebra over $\k$. For a $\pa$-closed element $a\in A$ we denote
by $[a]$ the corresponding cohomology class in $H^*(A):=H^*(A,\pa)$. 
Also for a homogeneous element
$a\in A$ we set
$$\ov{a}=(-1)^{1+\deg(a)}a.$$
Suppose that we have a collection of homogeneous elements 
$a_\bullet=(a_{ij})$, where $0\le i<j\le n$, $(i,j)\neq (0,n)$,
satisfying the equations
\begin{equation}\label{MC-eq}
\pa(a_{ij})=\sum_{i<k<j}\ov{a}_{ik}a_{kj}
\end{equation}
for all $0\le i<j\le n$, $(i,j)\neq (0,n)$ 
(in particular, the elements $a_{i,i+1}$ are $\pa$-closed). Then it is easy to check that
$$\mu(a_\bullet):=\sum_{0<k<n}\ov{a}_{0k}a_{kn}$$
is also $\pa$-closed. 
For given (homogeneous) cohomology classes $h_1,\ldots,h_n\in H^*(A)$,
one defines the $n$th Massey product
$$\lan h_1,\ldots,h_n\ran_{dg}\sub H^*(A)$$
as the subset formed by the classes $[\mu(a_\bullet)]$ as $a_\bullet=(a_{ij})$ runs through all collections
as above with $[a_{i-1,i}]=h_i$, $i=1,\ldots,g$ (see \cite{Kraines}, \cite{May}, \cite{LPWZ};
we follow the sign convention of \cite{LPWZ}).
We call a collection $a_\bullet$ as above a {\it defining system} for $\lan h_1,\ldots,h_n\ran_{dg}$.
We say that the Massey product $\lan h_1,\ldots,h_n\ran_{dg}$ is defined if this subset is nonempty,
i.e., there exists a defining system for $\lan h_1,\ldots,h_n\ran_{dg}$.
For example, the double Massey product is always defined and is given by the usual product, 
up to a sign. The triple Massey product $\lan h_1,h_2,h_3\ran_{dg}$
is defined if and and only if the double products $h_1h_2$ and $h_2h_3$ vanish in $H^*(A)$.

Now let $\CC$ be a dg-category over $\k$, and let $H^*(\CC)$ be the corresponding graded
category obtained by passing to cohomology on morphisms. 
Let $X_0,\ldots,X_n$ be objects of $\CC$.
The equations \eqref{MC-eq} make sense for a collection of (homogeneous) morphisms 
$a_{ij}\in\Hom_\CC^*(X_j,X_i)$. Thus, similarly to the case of a dg-algebra, one defines
the Massey product of a collection $(h_1,\ldots,h_n)$
of homogeneous morphisms in $H^*(\CC)$, where $h_i\in H^*\Hom_\CC(X_i,X_{i-1})$.

Recall that the homological perturbation theory provides a minimal $A_{\infty}$-structure on
$H^*(\CC)$ (see, e.g., \cite{Merkulov}). 
We will use the following important connection between the dg Massey products
and this $A_\infty$-structure.

\begin{thm}\label{LPWZ-thm} (\cite[Thm.\ 3.1]{LPWZ})
Consider a minimal  $A_{\infty}$-structure on
$H^*(\CC)$ obtained by homological perturbation theory. Let
$(h_i\in H^*\Hom_\CC(X_i,X_{i-1}))$, $i=1,\ldots,n$, be homogeneous elements, such that the
Massey product $\lan h_1,\ldots,h_n\ran_{dg}$ is defined. Then
$$(-1)^bm_n(h_1,\ldots,h_n)\in\lan h_1,\ldots,h_n\ran_{dg}, \ \text{ where}$$
$$b=1+\deg h_{n-1} + \deg h_{n-3} + \deg h_{n-5}+\ldots.$$
\end{thm}

Recall that for a dg-category $\CC$ one also has the dg-category of (left) dg $\CC$-modules
$\CC-\mod$, i.e., of dg-functors from $\CC$ to the dg-category of $\k$-complexes.
For a pair $M_1,M_2$ of dg $\CC$-modules the obvious identification gives isomorphisms
of complexes
\begin{equation}\label{Hom-shift-eq}
\Hom_\CC^\bullet(M_1,M_2[1])\simeq\Hom_\CC^\bullet(M_1,M_2)[1],\ \
\Hom_\CC^\bullet(M_1[-1],M_2)\simeq\Hom_\CC^\bullet(M_1,M_2)[1]',
\end{equation}
where for a complex $(C,d)$ we denote by $C[1]'$ the complex $(C[1],d)$ 
(recall that the usual shifted complex $C[1]$ has the differential $-d$).
We need to investigate the effect of shifts on the Massey products.

\begin{lem}\label{dg-Massey-shift-lem} 
Let $M_0,M_1,\ldots,M_n$ be dg $\CC$-modules, and let
$\wt{M}_i=M_i[k_i]$ for some $k_i\in\Z$, $i=0,\ldots,n$.
Let $h_i\in H^*\Hom_\CC(M_i,M_{i-1})$ be homogeneous elements, such that the
Massey product $\lan h_1,\ldots,h_n\ran_{dg}$ is defined.
Let $\wt{h}_i\in H^*\Hom_\CC(\wt{M}_i,\wt{M}_{i-1})$ be an element corresponding to $h_i$
under the obvious identification between the relevant spaces. Then one has
$$\lan\wt{h}_1,\ldots\wt{h}_n\ran_{dg}=(-1)^{k_1+\ldots+k_{n-1}}
\lan h_1,\ldots,h_n\ran\sub H^*\Hom_\CC(\wt{M}_n,\wt{M}_0)\simeq H^*\Hom_\CC(M_n,M_0).$$
\end{lem}

\Pf . It is enough to consider the case when $k_i=0$ for all $i\neq i_0$ and $k_{i_0}=1$. 
Let $a_\bullet=(a_{ij})$ be a defining system for $\lan h_1,\ldots,h_n\ran_{dg}$, so that
$[a_{i,i-1}]=h_i$. Set
$$\wt{a}_{ij}=\begin{cases} -a_{ij}, & i<i_0<j,\\ a_{ij}, & \text{otherwise.}\end{cases}$$
Using isomorphisms \eqref{Hom-shift-eq}, one can easily check that
$\wt{a}_\bullet=(\wt{a}_{ij})$ is a defining system for $\lan \wt{h}_1,\ldots,\wt{h}_n\ran_{dg}$.
Note that in the case $i_0=0$ or $i_0=n$ we have $\wt{a}_{ij}=a_{ij}$, so the two
Massey products are the same. On the other hand, if $0<i_0<n$ then
$$\mu(\wt{a}_\bullet)=-\mu(a_\bullet),$$
which implies the result.
\ed

Using Lemma \ref{dg-Massey-shift-lem} we can reduce the study of dg Massey products for the category of
dg-modules to the case when all $a_{i,i-1}$ have degree $1$ (and hence all $a_{ij}$ in a defining system have degree $1$). This will allow us to relate defining systems for Massey products
to twisted complexes.

\subsection{Convolutions in dg and triangulated categories}

Let $\CC$ be a dg-category.
A {\it twisted complex}\footnote{Our convention on the numbering of $M_i$ and degrees of $a_{ij}$ differs from that of
\cite{BK}.} 
over $\CC-\mod$ is a collection $\bM=(M_i,a_{ij})$, where
$M_i$, $i\in\Z$, are dg $\CC$-modules with $M_i=0$ for $i\gg 0$, and
$a_{ij}:M_j\to M_i$, $i<j$, are morphisms of degree $1$ satisfying
\begin{equation}\label{twisted-deg-1-eq}
\pa(a_{ij})+\sum_{i<k<j} a_{ij}a_{jk}=0
\end{equation}
for all $i<j$. 
We define the convolution of $\bM=(M_i,a_{ij})$ as the following $\CC$-module:
\begin{equation}\label{dg-convolution}
\conv(\bM):=(\bigoplus_i M_i, \pa+A),
\end{equation}
where $\pa$ is the usual differential on $\bigoplus_i M_i$ and 
$A=(a_{ij})$ is the upper-triangular endomorphism of $\bigoplus_i M_i$ with
components $a_{ij}$.
It is easy to check that $(\pa+A)^2=0$, so \eqref{dg-convolution} is indeed a dg $\CC$-module.

Let $f:M_1\to M_0$ be a closed morphism of degree $1$ in $\CC-\mod$. 
We can view this morphism as a twisted complex. Note that its convolution has form
$$\Cone_{dg}(f):=\conv(f)=(M_1\oplus M_2,\pa_{M_1\oplus M_2}+f).$$
We can also view $f$ as a closed morphism $\wt{f}:M_1[-1]\to M_0$ of degree $0$,
and $\Cone_{dg}(f)$ can be identified with the standard cone of $\wt{f}$.

The convolution of an arbitrary twisted complex can be obtained by iterating the
cone operation.

\begin{lem}\label{cone-lem}
Let $\bM=(M_i,a_{ij})$ be a twisted complex over $\CC-\mod$, where $M_i\neq 0$ only for
$i=0,\ldots,n$. Let us consider the truncated twisted complex $\tau_{[1,n]}\bM$ obtained
by considering only $M_1,\ldots,M_n$ and $a_{ij}$ with $i\ge 1$. Then $(a_{0i})$ are
components of a closed morphism of degree $1$
$$(a_{0\bullet}):\conv(\tau_{[1,n]}\bM)\to M_0,$$
and there is a natural isomorphism of dg $\CC$-modules
$$\conv(\bM)\simeq\Cone_{dg}(a_{0\bullet}).$$
\end{lem}

The proof is straightforward and is left to the reader.

Using the above lemma we can connect convolutions of twisted complexes over $\CC-\mod$ with convolutions in the triangulated category $H^0(\CC-\mod)$.
Recall (see \cite[Exer.\ IV.2.1]{GM}) that a {\it convolution} of a complex 
\begin{equation}\label{convolution-complex-eq}
X_n\rTo{d_n} X_{n-1}\to\ldots \to X_1 \rTo{d_1}X_0
\end{equation}
in a triangulated category $\TT$ (so that $d_i\circ d_{i+1}=0$), is an object $T\in\TT$ equipped
with morphisms $\a:T\to X_n[n]$, $\b:X_0\to T$, such that there exists a diagram
(called a {\it left Postnikov diagram})
\begin{diagram}
       &                     &X_{n-1}&\rTo(2,0)^{d_{n-1}}&&            &X_{n-2}           &\rTo(2,0)^{d_{n-2}}\ldots  &&X_0\\
       &\ruTo^{d_n}&              &\rdTo^{}             &     &\ruTo^{}&                        &      &\ruTo{}             &        &\rdTo^{\b}\\
X_n=C_n&\lTo(2,0)^{[1]} &&                     &C_{n-1}&\lTo(2,0)^{[1]}\ldots &&C_1&\lTo(2,0)^{[1]}&       &            & C_0=T
\end{diagram}
in which all the triangles $(X_i,C_i,C_{i+1})$ are distinguished, so that
$\a:T\to X_n[n]$ is the composition of the arrows in the lower row.

\begin{lem}\label{convolution-shift-lem} 
Let $(T,\a,\b)$ be a convolution of a complex \eqref{convolution-complex-eq}. Then
$(T[1],(-1)^n\a,\b)$ is a convolution of the shifted complex
$$X_n[1]\rTo{d_n} X_{n-1}[1]\to\ldots\to X_1[1]\rTo{d_1}X_0[1].$$
\end{lem}

\Pf . This can be easily checked by induction on $n$ using the fact that if
$$X\rTo{f} Y\rTo{g} Z\rTo{h} X[1]$$
is a distinguished triangle then the triangle
$$X[1]\rTo{f} Y[1]\rTo{g} Z[1]\rTo{-h} X[2]$$
is also distinguished.
\ed

\begin{lem}\label{convolution-lem}
Let $\bM=(M_i,a_{ij})$ be a twisted complex over $\CC-\mod$, where $M_i\neq 0$ only for
$i=0,\ldots,n$. Consider the complex
\begin{equation}\label{modules-complex-eq}
M_n[-n]\rTo{d_n} M_{n-1}[-n+1]\to\ldots\to M_1[-1]\rTo{d_1} M_0
\end{equation}
in the triangulated category $H^0(\CC-\mod)$, where $d_i=[a_{i-1,i}]$.
Let $\pi:\conv(\bM)\to M_n=(M_n[-n])[n]$ and $\iota:M_0\to\conv(\bM)$ be the natural maps
(given by the projection and the inclusion, respectively). 
Then the dg $\CC$-module $\conv(\bM)$, together with the maps $\a=(-1)^{{n\choose 2}}\pi$ and $\b=\iota$, is a convolution of
the complex \eqref{modules-complex-eq} in $H^0(\CC-\mod)$.
\end{lem}

\Pf . We can proceed by induction on $n$ (the case $n=1$ was discussed before). 
Let $T=\conv(\tau_{[1,n]}\bM)$. By the induction assumption,
$(T, (-1)^{{n-1\choose 2}}\pi':T\to M_n, \iota':M_1\to T)$ is a convolution of the complex 
$$M_n[-n+1]\rTo{d_n} M_{n-1}[-n+2]\to\ldots\rTo{d_2} M_1.$$
Hence, by Lemma \ref{convolution-shift-lem}, $(T[-1],(-1)^{n\choose 2}\pi',\iota')$ is a convolution of
$$M_n[-n]\rTo{d_n} M_{n-1}[-n+1]\to\ldots\rTo{d_2}M_1[-1].$$
On the other hand, by Lemma \ref{cone-lem},
$$\conv(\bM)=\Cone(\conv(\tau_{[1,n]}\bM)[-1]\rTo{a_{0\bullet}} M_0),$$
and the assertion follows.
\ed

\subsection{Massey products for triangulated categories}\label{Mas-prod-sec}

Let us recall the definition of the Massey products for triangulated categories, sometimes called 
Toda brackets (see \cite{Cohen}, \cite[Exer.\ IV.2.3]{GM}).
Let
\begin{equation}\label{triang-complex}
X_n\rTo{d_n} X_{n-1}\rTo{}\ldots \rTo{} X_1\rTo{d_1} X_0
\end{equation}
be a complex in a triangulated category (so that $d_i\circ d_{i+1}=0$).
One defines $\lan d_1,\ldots,d_n\ran\sub\Hom^{2-n}(X_n,X_0)$ as the set of all $p\circ q$,
where we take a convolution $(T, \a:T\to X_{n-1}[n-2],\b:X_1\to T)$ of the complex $X_{n-1}\to\ldots\to X_1$ (if it exists), and 
pick morphisms $p:T\to X_0$ and $q:X_n\to T[2-n]$, such that $d_n$ is equal to the composition
$X_n\rTo{q} T[2-n]\rTo{\a} X_{n-1}$ and $d_1$ is equal to the composition
$X_1\rTo{\b} T\rTo{p} X_0$.


It is well known that $0\in\lan d_1,\ldots,d_n\ran$ is exactly the condition for the existence of a convolution of the complex \eqref{triang-complex}. 
Also, for
$\lan d_1,\ldots,d_n\ran$ to be nonempty it is necessary that
$0\in\lan d_1,\ldots,d_{n-1}\ran$ and $0\in\lan d_2,\ldots,d_n\ran$ (and hence, the same
is true for any proper substring in $d_1,\ldots,d_n$).

It is easy to see that the triple product $\lan d_1,d_2,d_3\ran$ is nonempty provided
$d_2\circ d_3=d_1\circ d_2=0$, and is a coset for the subgroup
$$d_1\circ \Hom^{-1}(X_3,X_1)+\Hom^{-1}(X_0,X_2)\circ d_3\sub \Hom^{-1}(X_3,X_0).$$
Note that a similar result holds for triple dg Massey products.
The case of higher Massey products is more complicated. We will only consider a certain particular
situation for the quadruple products in the triangulated and dg-settings (see Lemma \ref{4Massey-amb-lem} below).
 
The following relation between the Massey products in
dg-categories and triangulated categories is well known to the experts and its various
versions have appeared in the literature (see \cite[Prop.\ 6.5]{Shipley} and \cite[Sec.\ 5.A]{BK}, which refers to the dissertation by Kapranov \cite{Kap}).

\begin{prop}\label{Massey-prop} 
Let $\CC$ be a dg-category, $M_0,\ldots,M_n$ a collection of dg $\CC$-modules, and 
$d_i\in H^{k_i}\Hom_\CC(M_i,M_{i-1})$, $i=1,\ldots,n$, a collection of maps in $H^*(\CC-\mod)$,
such that the dg Massey product $\lan d_1,\ldots,d_n\ran_{dg}$ is defined.
Then we have 
$$(-1)^{\sum_{i=1}^{n-1}(n-i)k_i}\lan d_1,\ldots,d_n\ran_{dg}\sub \lan d_1,\ldots,d_n\ran,$$
where on the right we consider the Massey product for the complex 
\begin{equation}\label{shifted-Massey-complex}
M_n[-k_1-\ldots-k_n]\rTo{d_n} M_{n-1}[-k_1-\ldots-k_{n-1}]\to\ldots\to M_1[-k_1]\rTo{d_1} M_0
\end{equation}
in the triangulated category $H^0(\CC-\mod)$.
\end{prop}

\Pf . By Lemma \ref{dg-Massey-shift-lem}, it is enough to consider the case when all $k_i=1$. In this
case we have to prove that
$$(-1)^{{n-1\choose 2}}\lan d_1,\ldots,d_n\ran_{dg}\sub \lan d_1,\ldots,d_n\ran,$$
where on the right we consider the Massey product for the complex
$$M_n[-n]\rTo{d_n} M_{n-1}[-n+1]\to\ldots\rTo{d_1} M_0.$$
Let $a_\bullet=(a_{ij})$, $a_{ij}\in\Hom^1_\CC(M_j,M_i)$ 
be a defining system for $\lan d_1,\ldots,d_n\ran_{dg}$,
so that $[a_{i-1,i}]=d_i$ and \eqref{MC-eq} is satisfied with
$\ov{a}_{ik}=a_{ik}$.
Considering the restricted system $(a_{ij}\ |\ 1\le i,j\le n-1)$
we obtain a twisted complex 
$$\bM=(M_{n-1},\ldots,M_1, (-a_{ij})_{1\le i<j\le n-1}).$$
Let $T=\conv(\bM)$ be the convolution of $\bM$. By Lemma \ref{convolution-lem}, 
$(T, (-1)^{{n-2\choose 2}}\pi:T\to M_{n-1},\iota:M_0\to T)$ is a convolution of the complex 
$$M_{n-1}[-n+2]\rTo{d_{n-1}}\ldots\rTo{d_2}M_1$$ 
in the triangulated category $H^0(\CC-\mod)$.
Hence, by Lemma \ref{convolution-shift-lem}, $(T[-1],(-1)^{{n-1\choose 2}}\pi,\iota)$ is a convolution of
$$M_{n-1}[-n+1]\rTo{d_{n-1}}\ldots\rTo{d_2}M_1[-1].$$ 
Furthermore, we have closed morphisms of degree $1$ of dg $\CC$-modules
$$\wt{q}=(a_{\bullet n}):M_n\to T, \ \ p=(a_{0\bullet}):T\to M_0$$
(this follows from \eqref{MC-eq} for $j=n$ and $i=0$, respectively), such that
$\pi\circ\wt{q}=a_{n-1,n}$ and $p\circ\iota=a_{01}$.
Thus, the morphisms in $H^0(\CC-\mod)$,
$$q=(-1)^{{n-1\choose 2}}\wt{q}:M_n[-n]\to T[-n+1]=T[-1][2-n] \ \text{ and } p:T[-1]\to M_0,$$ 
satisfy the conditions in the definition of the Massey product $\lan d_1,\ldots,d_n\ran$.
Since
$p\circ q\in\lan d_1,\ldots,d_n\ran$ is represented by $(-1)^{{n-1\choose 2}}\mu(a_\bullet)$, the assertion follows.
\ed

On the other hand, by Theorem \ref{LPWZ-thm}, 
the Massey product $\lan d_1,\ldots,d_n\ran_{dg}$ always contains
$\pm m_n(d_1,\ldots,d_n)$, where $(m_\bullet)$ is a minimal $A_{\infty}$-structure
on $H^*(\CC)$ obtained by homological perturbation theory. This leads to
the following result that
will allow us to compute the Hochschild cohomology class of $m_3$ and, in a special situation,
of $m_4$, via the Massey products.

\begin{cor}\label{mn-Massey-cor}
In the situation of Proposition \ref{Massey-prop} consider a minimal  $A_{\infty}$-structure on
$H^*(\CC-\mod)$ obtained by the homological perturbation theory. 
Assume that the Massey product
$\lan d_1,\ldots, d_n\ran_{dg}$ is defined.
Then 
$$(-1)^{b+\sum_{i=1}^{n-1}(n-i)k_i}m_n(d_{1},\ldots,d_{n})\in \lan d_{1},\ldots,d_{n}\ran$$
with $b=1+k_{n-1}+k_{n-3}+k_{n-5}+\ldots$,
where on the right we consider the Massey product for the complex \eqref{shifted-Massey-complex} in $H^0(\CC-\mod)$.
\end{cor}

\begin{rem}\label{enhanced-rem}
Any enhanced triangulated category in the sense of \cite{BK} 
can be realized as a full subcategory in $H^0(\CC-\mod)$ for the corresponding dg-category $\CC$
(e.g., this follows from \cite[Prop.\ 1.3, Prop.\ 3.2]{BK}).
Therefore, we can apply Corollary \ref{mn-Massey-cor} to compare the Massey products in an enhanced triangulated
category with the higher products obtained by the homological perturbation theory.
\end{rem}

Later we will need the following result. 
 
\begin{lem}\label{4Massey-amb-lem} 
(i) Suppose we have a complex \eqref{triang-complex} in a triangulated category, where $n=4$. Assume that $0\in\lan d_1,d_2,d_3\ran$, $0\in\lan d_2,d_3,d_4\ran$,
and a left Postnikov system 
for the complex
$$X_3\rTo{d_3} X_2\rTo{d_2} X_1$$ 
is unique up to an isomorphism (identical on $X_i$).
Then the Massey product $\lan d_1,d_2,d_3,d_4 \ran$ is nonempty.
Also, for any $\mu,\mu'\in \lan d_1,d_2,d_3,d_4 \ran$ one has
\begin{equation}\label{4Mas-diff-eq}
\mu-\mu'\in \lan \Hom^{-1}(X_2,X_0),d_3,d_4\ran+d_1\circ\Hom^{-2}(X_4,X_1).
\end{equation}

\noindent
(ii) Let $\CC$ be a dg-category, and let $d_1,d_2,d_3,d_4$ be
a sequence of composable arrows in $H^0(\CC-\mod)$ satisfying the assumptions of (i). Then
the dg Massey product $\lan d_1,d_2,d_3,d_4\ran_{dg}$ is defined.
\end{lem}

\Pf . (i) By definition, $\lan d_1,d_2,d_3,d_4\ran$ consists of $p\circ q$, where $p$ and $q$
come from a diagram
\begin{equation}\label{4Mas-diagram}
\begin{diagram}
X_4&\rTo{d_4}&X_3&\rTo{d_3}&X_2&\rTo{d_2}&X_1&\rTo^{d_1}&X_0\\
&\rdTo(4,4)^{q}_{[-2]}\rdTo(4,2)^{[-1]}_{t}&&\luTo^{[1]}_{\pi}&\dTo{\iota}&\ruTo^{\wt{d}_2}\ldTo(2,4){\de}&&
\ruTo(4,4)_{p}\\
&& &&P\\
&&&&\uTo_{\eps}^{[1]}\\
&&&& T
\end{diagram}
\end{equation}
in which the triangles $(X_3,X_2,P)$ and $(P,X_1,T)$ are distinguished (so in the middle we have a left Postnikov system for $X_3\to X_2\to X_1$) and all other
triangles are commutative. To show the existence of such a diagram we observe first that
we can always construct a left Postnikov system in the middle and a morphism $t$ such that
$\pi\circ t=d_4$ (since $d_3\circ d_4=0$). We have
$$\wt{d}_2\circ t\in\lan d_2,d_3,d_4\ran.$$
Hence, the assumption $0\in\lan d_2,d_3,d_4\ran$ implies 
$$\wt{d}_2\circ t=d_2\circ f+g\circ d_4$$
for some $f\in\Hom^{-1}(X_4,X_2)$ and $g\in\Hom^{-1}(X_3,X_1)$. 
Thus, changing $\wt{d}_2$ to $\wt{d}_2-g\circ\pi$ and $t$ to $t-\iota\circ f$, we can
achieve that 
\begin{equation}\label{t-choice-eq}
\pi\circ t=d_4, \ \wt{d}_2\circ t=0.
\end{equation}
By the uniqueness of the left Postnikov diagram in the middle, in fact, the needed $t$ exists for
any choice of $\wt{d}_2$ (since two such choices differ by an automorphism of $P$ compatible with $\pi$ and
$\iota$). Once we have $t$ satisfying \eqref{t-choice-eq}, we can find $q$ such that $\eps\circ q=t$.
On the other hand, the morphism $p$ in the diagram exists provided
$d_1\circ\wt{d}_2=0$. It is easy to see that
$$d_1\circ\wt{d}_2=\mu\circ\pi$$
for some $\mu\in\lan d_1,d_2,d_3\ran$. Hence,
$$\mu=d_1\circ g'+h\circ d_3$$
for some $g'\in\Hom^{-1}(X_3,X_1)$ and $h\in\Hom^{-1}(X_2,X_0)$. This implies that
$$d_1\circ\wt{d}_2=d_1\circ g'\circ\pi,$$
so, changing $\wt{d}_2$ to $\wt{d}_2-g'\circ\pi$, we will have
$d_1\circ\wt{d}_2=0$, which will give the morphism $p$.

It remains to establish \eqref{4Mas-diff-eq}.
By assumption, up to an isomorphism, two diagrams \eqref{4Mas-diagram} differ
only by a choice of the maps $(t,p,q)$. Given another diagram with maps $(t',p',q')$ we can
write
$$p'\circ q'-p\circ q=(p'-p)\circ q'+p\circ (q'-q).$$
Now we have $q'-q=\de\circ x$ for some $x\in\Hom^{-2}(X_4,X_1)$ and
$p'-p=y\circ\eps$ for some $y\in\Hom^{-1}(P,X_0)$. Therefore,
$$p'\circ q'-p\circ q=y\circ t+d_1\circ x.$$
It remains to observe that $y\circ t\in\lan y\circ\iota,d_3,d_4\ran$.

\noindent
(ii) Let $a_{i-1,i}\in\Hom_{\CC}^0(X_i,X_{i-1})$ be representatives of $d_i$ for $i=1,2,3,4$.
By assumption, there exist elements $a_{24}\in\Hom_\CC^{-1}(X_4,X_2)$ and $a_{13}\in\Hom_\CC^{-1}(X_3,X_1)$ such that 
$$\pa(a_{24})=-a_{23}a_{34}, \ \pa(a_{13})=-a_{12}a_{23}.$$
By Proposition \ref{Massey-prop}, we have $0\in\lan d_2,d_3,d_4\ran_{dg}$. Hence, 
$$a_{13}a_{34}-a_{12}a_{24}=xa_{34}+a_{12}y+\pa(a_{14})$$
for some $a_{14}\in \Hom^{-2}_\CC(X_4,X_1)$, $x\in \Hom^{-1}_\CC(X_3,X_1)$ and 
$y\in \Hom^{-1}_\CC(X_4,X_2)$, such that $\pa(x)=0$, $\pa(y)=0$.
Hence, changing $a_{13}$ to $a_{13}-x$ and $a_{24}$ to $a_{24}+y$, we can achieve
that 
$$a_{13}a_{34}-a_{12}a_{24}=\pa(a_{14}).$$
Similarly, from the condition $0\in\lan d_1,d_2,d_3\ran$ we obtain that for some
$a'_{13}$, $a_{02}$ and $a_{03}$ one has
$$\pa(a'_{13})=-a_{12}a_{23},\ \pa(a_{02})=-a_{01}a_{12},$$
$$a_{02}a_{23}-a_{01}a'_{13}=\pa(a_{03}).$$
Now we need to use our assumption on the uniqueness of a Postnikov system, up to
an isomorphism, to find a relation between $a'_{13}$ and $a_{13}$.
Let 
$$P=\Cone_{dg}(a_{23})\in\CC-\mod$$ 
be the cone of $a_{23}$, viewed as a closed morphism of degree $1$ from $X_3[1]$ to $X_2$, so that
we have a triangle of closed morphisms
$$X_3\rTo{a_{23}} X_2\rTo{\iota} P\rTo{\pi} X_3[1]$$
that becomes distinguished in the triangulated category $H^0(\CC-\mod)$. Our assumpion on the uniqueness of a Postnikov system
means that there exists a unique
morphism $\wt{d}_2\in H^0\Hom_{\CC}(P,X_1)$ such that $\iota\circ\wt{d}_2=d_2$ in 
$H^0\Hom_{\CC}(X_2,X_1)$, 
up to an automorphism of $P$ in $H^0(\CC-\mod)$, compatible with the cone structure of $P$.
We have two such morphisms $\wt{d}_2$, namely
$$\wt{d}_2=(-a_{12},a_{13}) \mod \im(\pa) \ \text{ and } \wt{d}'_2=(-a_{12},a'_{13}) \mod \im(\pa).$$
Therefore, we have 
\begin{equation}\label{F-dd-eq}
\wt{d}'_2=\wt{d}_2\circ F
\end{equation}
for some automorphism $F:P\to P$ in $H^0(\CC-\mod)$, 
compatible with the cone structure of $P$. Any such automorphism has form 
$$F=\id_P-\iota f \pi \mod \im(\pa)$$
for some closed element $f\in \Hom_\CC^{-1}(X_3,X_2)$.
Hence, the condition \eqref{F-dd-eq} gives
$$a'_{13}=a_{13}+a_{12} f+g a_{23}+\pa(h),$$
where $g\in\Hom^{-1}_\CC(X_2,X_1)$, $\pa(g)=0$, and $h\in \Hom^{-2}_\CC(X_3,X_1)$
(the term $g a_{23}$ comes from the form of the differential on $\Hom_\CC(P,X_1)$).
Now setting 
$$a'_{02}=a_{02}-a_{01}g,\ \ a'_{03}=a_{03}-a_{02}f+a_{01}h$$
we obtain that
$$(a_{01},a_{12},a_{23},a_{34},a'_{02},a_{13},a_{24},a'_{03},a_{14})$$
is a defining system for $\lan d_1,d_2,d_3,d_4\ran$.
\ed

\subsection{Some triple Massey products on curves}\label{3Mas-sec}

Let $C$ be a curve and $p\in C$ a smooth point.
Let us denote by $\xi_p$ the image of a generator of $H^0(\OO(p)/\OO)$ under 
the connecting homomorphism
$H^0(\OO(p)/\OO)\to H^1(\OO)$. 
We would like to study the map
\begin{equation}\label{Mas-prod-eq}
\Ext^1(\OO_p,\OO)\ot\Ext^1(\OO_p,\OO_p)\ot\Hom(\OO,\OO_p)\to H^1(C,\OO)/\lan\xi_p\ran
\end{equation}
given by the triple Massey product in $D^b(C)$ of the type
$$\OO[-2]\to \OO_p[-2] \to \OO_p[-1] \to \OO.$$
Note that such a Massey product is always nonempty since $\Ext^1(\OO,\OO_p)=\Ext^2(\OO_p,\OO)=0$
and the ambiguity is exactly $\lan\xi_p\ran\sub H^1(C,\OO)$, which is equal to the image of the composition map
$$\Ext^1(\OO_p,\OO)\ot\Hom(\OO,\OO_p)\to H^1(C,\OO).$$
It is also compatible with the map $-m_3$, obtained by the homological perturbation theory (see Corollary \ref{mn-Massey-cor} and
Remark \ref{enhanced-rem}).

Note that we have canonical bases in spaces
$\Hom(\OO,\OO_p)$, $\Ext^1(\OO_p,\OO_p)\ot T^*$ and $\Ext^1(\OO_p,\OO)\ot T^*$,
where $T=T_p$ is the tangent line to $C$ at $p$. 
By the definition of the Massey product, we have to consider the 
canonical extension
$$0\to T^*\ot_{\k}\OO_p\rTo{i} \OO_{2p}\rTo{\pi} \OO_p\to 0$$
inducing a generator of $\Ext^1(\OO_p,\OO_p)$.
Then we should consider
liftings $c:\OO\to \OO_{2p}$ and $d:\OO_{2p}\to (T^*)^{\ot 2}\ot_{\k}\OO[1]$ such that
$\pi\circ c=1:\OO\to\OO_p$ and $i\circ d: T^*\ot_{\k}\OO_p\to (T^*)^{\ot 2}\ot_{\k}\OO[1]$
is the canonical element represented by the extension
$$0\to (T^*)^{\ot 2}\ot_{\k}\OO\to (T^*)^{\ot 2}\ot_{\k}\ot\OO(p)\to T^*\ot_{\k}\OO_p\to 0.$$
Thus, we can take $c=1\in H^0(C,\OO_{2p})$ and $d$ to be the class of the extension
\begin{equation}\label{ext-d-eq}
0\to (T^*)^{\ot 2}\ot_{\k}\OO\to (T^*)^{\ot 2}\ot_{\k}\OO(2p)\to \OO_{2p}\to 0.
\end{equation}
Our Massey product is the coset of the
composition $d\circ c:\OO\to (T^*)^{\ot 2}\ot_{\k}\OO[1]$
in $(T^*)^{\ot 2}\ot_{\k} H^1(\OO)/\lan\xi_p\ran$.
In other words, this is the image of $1\in H^0(\OO_{2p})$ under the
boundary homomorphism
\begin{equation}\label{bound-hom-eq}
\de_{2p}: H^0(\OO_{2p})\to (T^*)^{\ot 2}\ot_{\k} H^1(\OO)
\end{equation}
associated with the extension \eqref{ext-d-eq}, viewed
modulo $\lan\xi_p\ran$. Since the latter subspace
is the image under $\de_{2p}$
of the subspace $H^0(T^*\ot_{\k}\OO_p)\sub H^0(\OO_{2p})$, we obtain
that our Massey product is zero if and only if 
$$H^0(\OO_{2p})=\ker(\de_{2p})+H^0(T^*\ot_{\k}\OO_p).$$
Since $\ker(\de_{2p})$ is the image of the homomorphism
$$(T^*)^{\ot 2}\ot_{\k} H^0(\OO(2p))\to H^0(\OO_{2p}),$$
the Massey product vanishes if and only if the composed map
$$(T^*)^{\ot 2}\ot_{\k} H^0(\OO(2p))\to H^0(\OO_{2p})\to 
H^0(\OO_{2p})/H^0(T^*\ot_{\k}\OO_p)\simeq H^0(\OO_p)$$
is surjective. In other words, this is equivalent to surjectivity of the map
$$H^0(\OO(2p))\to H^0(\OO(2p)/\OO(p)),$$
or to the condition $H^0(C,\OO(2p))\not\sub H^0(C,\OO(p))$. Thus, we obtain the following result.

\begin{prop}\label{van-Mas-prop} 
Let $C$ be a curve, $p\in C$ a smooth point. 
The Massey product \eqref{Mas-prod-eq} vanishes if and only if
$H^0(C,\OO(2p))\not\sub H^0(C,\OO(p))$.
For example, if $C$ is smooth and projective of genus $g\ge 1$
then this happens if and only if
either $g=1$ or $C$ is hyperelliptic and $p$ is a Weierstrass point of $C$.
\end{prop}

Next, we are going to compute the Massey product \eqref{Mas-prod-eq} in terms
of additional data allowing us to represent classes in $H^1(\OO)$.
Namely, let $g$ be the arithmetic genus of $C$, and 
let us assume that $D$ is an effective divisor of degree $g-1$ (supported on the smooth
part of $C$) such that
$h^0(D+p)=1$ and $p\not\in\supp(D)$. Then the boundary homomorphism
$$\de_{D+p}:H^0(\OO(D+p)/\OO)\to H^1(\OO)$$
associated with an exact sequence $0\to \OO\to \OO(D+p)\to \OO(D+p)/\OO\to 0$
is an isomorphism.
Consider also the similar boundary homomorphism
$$\de_{D+2p}:H^0(\OO(D+2p)/\OO)\to H^1(\OO),$$
so that $\de_{D+p}$ is the restriction of $\de_{D+2p}$ to
the subspace $H^0(\OO(D+p)/\OO)\sub H^0(\OO(D+2p)/\OO)$.
Note that the kernel of $\de_{D+2p}$ is the image of the natural embedding
$$H^0(\OO(D+2p))/H^0(\OO)\to H^0(\OO(D+2p)/\OO).$$
Thus, for $x\in H^0(\OO(D+2p)/\OO)$ we can write
$$\de_{D+2p}(x)=\de_{D+p}(y)$$ where $y\in H^0(\OO(D+p)/\OO)$ is such that
$x\equiv y+s \mod\OO$ for some global section
$s\in H^0(\OO(D+2p))$.

We want to compute the image of a generator of $H^0(\OO(2p)/\OO)$ under
the composition of the boundary homomorphism 
$$\de_{2p}:H^0(\OO(2p)/\OO)\to H^1(\OO),$$
with the projection $H^1(\OO)\to H^1(\OO)/\lan\xi_p\ran$.
Since $\de_{2p}$ is just the restriction of $\de_{D+2p}$, we can apply the above recipe
to $x\in H^0(\OO(2p)/\OO)\sub H^0(\OO(D+2p)/\OO)$.
Note that $\lan\xi_p\ran=\de_{D+p}(H^0(\OO(p)/\OO))$, so we need to
find $y\in H^0(\OO(D+p)/\OO)$ and $s\in H^0(\OO(D+2p))$ such that
$$x\equiv y+s \mod\OO$$
and then view $y$ modulo $H^0(\OO(p)/\OO)$. In other words, we need to consider the
projection of $y$ to $H^0(\OO(D+p)/\OO(p))\simeq H^0(\OO(D)/\OO)$.
Since the polar part of $y$ near $\supp D$ is opposite to that of $s$, we obtain the following
formula for the Massey product \eqref{Mas-prod-eq}.

\begin{prop}\label{3-Mas-prod-prop} 
With the above choice of divisor $D$ let us consider the restriction maps
$$r_{D+2p,p}:
H^0(\OO(D+2p))/H^0(\OO)\rTo{\sim} H^0(\OO(D+2p)/\OO(D+p))\simeq T^{\ot 2} \ \ \text{and}$$
$$r_{D+2p,D}: H^0(\OO(D+2p))/H^0(\OO)\to H^0(\OO(D)/\OO).$$
Then the map \eqref{Mas-prod-eq} is equal to 
$$-\ov{\de_{D+p}}\circ r_{D+2p,D}\circ r_{D+2p,p}^{-1}: T^{\ot 2}\to H^1(\OO)/\lan\xi_p\ran,$$
where
$$\ov{\de_{D+p}}:H^0(\OO(D)/\OO)\rTo{\sim} H^0(\OO(D+p)/\OO(p))\to H^1(\OO)/\lan\xi_p\ran$$
is the isomorphism induced by $\de_{D+p}$.
\end{prop}

Assume now that we are in the situation of Section \ref{Egn-sec} with $n=g$, so we have
$g$ distinct smooth points $p_1,\ldots,p_g\in C$
such that $h^0(p_1+\ldots+p_g)=1$, and the corresponding classes
$X_i$, $i=1,\ldots,g$, form a basis in $H^1(C,\OO)$ 
(we use the notation from Section \ref{Egn-sec} for the basis elements in various $\Ext$-spaces).
Let $T_{p_i}$ denote the tangent line to $C$ at $p_i$. We have a natural isomorphism
$$T_{p_i}\simeq H^0(C,\OO(p_i)/\OO)\simeq\Ext^1(\OO_{p_i},\OO_{p_i}),$$
so we can think of $Y_i$ as a generator of $T_{p_i}$.
Let us set $D_i=\sum_{j\neq i} p_j$.

\begin{cor}\label{alpha-for-cor} 
The constants $\a_{ij}(m_3)$ associated with 
the natural $A_\infty$-structure on $E_{g,g}$ 
(see Proposition \ref{Hoch-prop}) can be computed as follows.
Pick an element $\wt{Y_i}\in H^0(\OO(2p_i+D_i))$
such that 
$$\wt{Y_i}\mod \OO(p_i+D_i)=(Y_i)^{\ot 2}\in H^0(\OO(2p_i)/\OO(p_i))\simeq T_{p_i}^{\ot 2}.$$
Then
$$\a_{ij}(m_3)\cdot Y_j=\wt{Y_i}\mod \OO(2p_i+\sum_{k\neq i,j}p_k)\in H^0(\OO(p_j)/\OO)
\simeq T_{p_j}.$$
\end{cor}

\Pf . This follows from the above computation of the Massey product $\lan B_i,Y_i,A_i\ran$ together with the compatibility
$$-m_3(B_i,Y_i,A_i)\in \lan B_i,Y_i,A_i\ran$$
obtained from Corollary \ref{mn-Massey-cor}.
\ed

\begin{rem} The above Corollary shows that the constants $\a_{ij}(m_3)$ are related
to a different kind of triple Massey product in $D^b(C)$ studied in \cite{P-Mas}. Namely, 
setting $D=\sum_{k=1}^g p_k$, we have 
\begin{equation}\label{two-kinds-Massey-eq}
\a_{ij}(m_3)\ot Y_j=\lan m_3(\OO(D),p_i,p_j), Y_i^{\ot 2}\ran,
\end{equation}
where $m_3(L,x,y)\in(\om_C\ot L^{-1})|_x\ot L|_y$ 
is the triple Massey product corresponding to the composable arrows
$$\OO_C\to\OO_x\rTo{[1]} L\to \OO_y$$
defined whenever $x\neq y$, $x$ is a base point of $\om_C\ot L^{-1}$ and $y$ is a base point of $L$
(see \cite[Sec.\ 1.1]{P-Mas}). Note that in our case
$$m_3(\OO(D),p_i,p_j)\in T^*_{p_i}\ot\OO(-D)|_{p_i}\ot\OO(D)|_{p_j}\simeq (T^*_{p_i})^{\ot 2}\ot 
T_{p_j}.$$
The identity \eqref{two-kinds-Massey-eq} also follows from the $A_{\infty}$-axioms associated with composable arrows 
$$\OO\to \OO_{p_l}\rTo{[1]} \OO(D)\to \OO_{p_k}\rTo{[1]}\OO$$
for $l=k$ and $l=i$.
Picking one more point $q\in C$ generically, we can write a formula for $m_3(\OO(D),p_i,p_j)$
in terms of theta-functions. Namely, first one easily checks that 
$m_3(\OO(D),p_i,p_j)=m_3(\OO(D-q),p_i,p_j)$.
Next, we represent $D-q$ as the sum of two divisors: 
$$D-q=\xi+(D_{ij}+q),$$ 
where $\xi=p_i+p_j-2q$ and $D_{ij}=\sum_{k\neq i,j}p_k$. 
We have a theta-function $\th_{D_{ij}+q}$ on the Jacobian
of $C$ associated with the degree-$(g-1)$ divisor $D_{ij}+q$. Now, by \cite[Lem.\ 2.2]{P-Mas},
we obtain
$$m_3(\OO(D),p_i,p_j)=m_3(\OO(D-q),p_i,p_j)=-\frac{\th_{D_{ij}+q}(2p_i-2q)\th'_{D_{ij}+q}(0)(p_i)}
{\th_{D_{ij}+q}(p_i-p_j)\th_{D_{ij}+q}(p_i+p_j-2q)},$$
where we view $\th'_{D_{ij}+q}(0)$ as a global $1$-form on $C$.
Note that the formula of \cite[Lem.\ 2.2]{P-Mas} is applicable since $p_i$ and $p_j$ are
not in the support of $D_{ij}+q$.
\end{rem}

\subsection{A quadruple Massey product}

Next, let us assume that $(C,p)$ is such that $H^0(C,\OO(2p))\not\sub H^0(C,\OO(p))$
(e.g., $C$ is a hyperelliptic smooth projective curve and $p$ is a Weierstrass point).

\begin{lem}\label{4Mas-def-lem} 
Under the above assumption the quadruple Massey product in $D^b(C)$ of the type
$$\OO[-3]\to\OO_p[-3]\to\OO_p[-2]\to\OO[-1]\rTo{\xi_p}\OO$$ 
gives rise to a well-defined map
\begin{equation}\label{4-Mas-prod-eq}
\lan\xi_p\ran\ot\Ext^1(\OO_p,\OO)\ot
\Ext^1(\OO_p,\OO_p)\ot\Hom(\OO,\OO_p)\to H^1(C,\OO)/\lan\xi_p\ran.
\end{equation}
The corresponding dg Massey product (coming from
some dg-enhancement of $D^b(C)$) is also defined.
\end{lem}

\Pf . We would like to apply Lemma \ref{4Massey-amb-lem} in our situation.
Note that the relevant triple Massey products contain zero, 
since $\Hom^2(\OO_p,\OO)=0$ and the triple Massey product \eqref{Mas-prod-eq} vanishes
by Proposition \ref{van-Mas-prop} (here we use our assumption on $(C,p)$).
Next, we need to check the uniqueness of a left Postnikov diagram (up to an isomorphism)
for the complex
$$\OO_p\rTo^{[1]}_{\psi}\OO_p\rTo^{[1]}_{\eta}\OO.$$
Since the distinguished triangle containing $\psi$ corresponds to a nontrivial extension
$$0\to\OO_p\rTo{\iota}\OO_{2p}\rTo{\pi}\OO_p\to 0,$$ 
it is enough to check that any diagram
\begin{diagram}
\OO_p&\rTo^{[1]}_{\psi}&\OO_p&\rTo^{[1]}_{\eta}&\OO&\\
&\luTo{\pi}&\dTo{\iota}&\ruTo_{[1]}^{\wt{\eta}}\\
&&\OO_{2p}
\end{diagram}
in which the left triangle is distinguished, is obtained from any other such diagram by an automorphism
of $\OO_{2p}$. Indeed, two choices of $\wt{\eta}$ differ by a morphism of the form
$f\circ\pi$, where $f\in\Hom^1(\OO_p,\OO)$. Thus, $f$ is a multiple of $\eta$: $f=c\cdot\eta$, and
so 
$$f\circ\pi=c(\eta\circ\pi)=\wt{\eta}\circ(c\cdot\iota\circ\pi).$$
Now consider the automorphism $\id+c(\iota\circ\pi)$ of $\OO_{2p}$
This automorphism is compatible with the extension structure and sends
$\wt{\eta}$ to $\wt{\eta}+f\circ\pi$, as required. 

It remains to check that the ambiguity for our Massey product is exactly 
$$\lan\xi_p\ran\sub H^1(C,\OO)=\Ext^1(\OO,\OO).$$ 
By Lemma \ref{4Massey-amb-lem}, we have to look at
the composition $\xi_p\circ\Hom(\OO,\OO)\sub \Ext^1(\OO,\OO)$ and at the
triple Massey product
$$\lan\Ext^1(\OO_p,\OO),\Ext^1(\OO_p,\OO_p),\Hom(\OO,\OO_p)\ran\sub 
\Ext^1(\OO,\OO).$$
Applying Proposition \ref{van-Mas-prop} once more, we see that the latter product is
$\lan\xi_p\ran$.

The last assertion follows from Lemma \ref{4Massey-amb-lem}(b).
\ed

By definition, the Massey product \eqref{4-Mas-prod-eq} is calculated as the composition
$\wt{\xi}\circ s$ in the diagram
\begin{equation}
\begin{diagram}
\OO&\rTo{}&\OO_p&\rTo^{[1]}_{\psi}&\OO_p&\rTo^{[1]}_{\eta}&\OO&\rTo^{[1]}_{\xi_p}&\OO\\
&\rdTo(4,4){s}\rdTo(4,2){}&&\luTo{}&\dTo{}&\ruTo_{[1]}^{\wt{\eta}}\ldTo(2,4){}&&
\ruTo(4,4)^{[1]}_{\wt{\xi}}\\
&& &&\OO_{2p}\\
&&&&\uTo{r}\\
&&&& \OO(2p)
\end{diagram}
\end{equation}
in which the triangle containing $\psi$ and
the triangle containing $\wt{\eta}$ and $r$ are distinguished and all the 
other triangles are commutative. Here $r$ is the composition of the natural
projection $\OO(2p)\to \OO(2p)/\OO$ and an isomorphism $\OO(2p)/\OO\simeq \OO_{2p}$.
In other words, we pick an element $\wt{\xi}\in H^1(\OO(-2p))$ such that
$t_*(\wt{\xi})=\xi_p$, where
$$t_*:H^1(\OO(-2p))\to H^1(\OO)$$
is the map induced by the canonical embedding $t:\OO(-2p)\to \OO$.
On the other hand, we choose a section $s:\OO\to\OO(2p)$ such that $r(s)=1\in H^0(\OO_{2p})$,
and apply the induced map $s_*:H^1(\OO(-2p))\to H^1(\OO)$ to $\wt{\xi}$.
One can check directly 
that the ambiguities in the choices of $\wt{\xi}$ and $s$ do not change the coset
of $s_*(\wt{\xi})$ in $H^1(\OO)/\lan\xi_p\ran$ (we also know this by Lemma \ref{4Mas-def-lem}). 
From this we obtain the following descriptions of the quadruple 
Massey product \eqref{4-Mas-prod-eq}
similar to those for the triple Massey product \eqref{Mas-prod-eq}.




\begin{prop}\label{4-Mas-prod-prop} 
Let $C$ be a curve, $p\in C$ a smooth point, such that 
$H^0(C,\OO(2p))\not\sub H^0(C,\OO(p))$. Let $T=T_p$ denote the tangent line to $C$ at $p$.

\noindent (a)
The Massey product \eqref{4-Mas-prod-eq} is given by the map
\begin{equation}\label{4-Mas-prod-for}
\phi: T^{\ot 2}\ot\lan\xi_p\ran\simeq T^{\ot 3}\simeq H^0(\OO(3p)/\OO(2p))\to
H^1(\OO)/\lan\xi_p\ran,
\end{equation}
where the last arrow is induced by the boundary homomorphism
$\de_{3p}:H^0(\OO(3p)/\OO)\to H^1(\OO)$.
The map $\phi$ vanishes if and only if $H^0(C,\OO(3p))\not\sub H^0(C,\OO(2p))$.

\noindent (b) Let $g$ be the arithmetic genus of $C$, and let 
$D$ be an effective divisor of degree $g-1$ (supported on the smooth part of 
$C$) such that $h^0(D+p)=1$ and
$p\not\in\supp(D)$. Then we have
$$\phi=-\ov{\de_{D+p}}\circ r_{D+3p,D}\circ r_{D+3p,p}^{-1},$$
where
$$r_{D+3p,p}: H^0(\OO(D+3p))/H^0(\OO(2p))\rTo{\sim} H^0(\OO(3p)/\OO(2p))\simeq T^{\ot 3},$$
$$r_{D+3p,D}: H^0(\OO(D+3p))/H^0(\OO(2p))\to H^0(\OO(D)/\OO)$$
are natural restriction maps
and $\ov{\de_{D+p}}$ is the isomorphism from Proposition \ref{3-Mas-prod-prop}.
\end{prop}

\Pf . (a) Pick $s\in H^0(\OO(2p))$ such that $\ov{s}=s\mod\OO(p)\neq 0$.
Then, as we have seen above,
$$\phi(\ov{s}\ot \xi_p)=s_*(t_*)^{-1}(\xi_p)$$
(the right-hand side is well defined in $H^1(\OO)/\lan\xi_p\ran$).
Recall that $\xi_p\in H^1(\OO)$ is the image of a generator $\psi_p$ of
$T\simeq H^0(\OO(p)/\OO)$ under the boundary map $H^0(\OO(p)/\OO)\to H^1(\OO)$.
Hence, if we pick a generator $\wt{\psi_p}\in H^0(\OO(p)/\OO(-2p))$, such that
$\wt{\psi_p}\equiv \psi_p\mod\OO$, then
$(t_*)^{-1}(\xi_p)\sub H^1(\OO(-2p))$ is represented modulo $\ker(t_*)$
by the image of $\wt{\psi_p}$
under the boundary homomorphism
$$H^0(\OO(p)/\OO(-2p))\to H^1(\OO(-2p)).$$
Now the morphism of exact sequences 
\begin{diagram}
0&\rTo{} &\OO(-2p)&\rTo{}&\OO(p)&\rTo{}&\OO(p)/\OO(-2p)&\rTo{}& 0\\
&&\dTo{s}&&\dTo{s}&&\dTo{s'}\\
0&\rTo{} &\OO&\rTo{}&\OO(3p)&\rTo{}&\OO(3p)/\OO&\rTo{}& 0
\end{diagram}
shows that $\phi(\ov{s}\ot \xi_p)$ is represented by the image of $s'(\wt{\psi_p})$
under the boundary map $H^0(\OO(3p)/\OO)\to H^1(\OO)$, which implies our first assertion.

The map $\phi$ vanishes if and only if the image of the boundary homomorphism
$\de_{3p}:H^0(\OO(3p)/\OO)\to H^1(\OO)$
is equal to $\lan \xi_p\ran$, which is the image of $\de_{2p}:H^0(\OO(2p)/\OO)\to H^1(\OO)$.
Since $H^0(\OO(2p)/\OO)$ has codimension $1$ in $H^0(\OO(3p)/\OO)$, this happens
exactly when $\ker(\de_{2p})$ has codimension $1$ in $\ker(\de_{3p})$.
But these kernels are $H^0(\OO(2p))/H^0(\OO)$ and $H^0(\OO(3p))/H^0(\OO)$,
respectively, hence the assertion.

\noindent
(b) Note that $\chi(D+p)=1$, so $H^1(\OO(D+p))=0$. It follows that $H^1(\OO(D+2p))=0$,
and hence, $h^0(D+2p)=\chi(D+2p)=2$. Also, $h^0(p)\le h^0(D+p)=1$, so $h^0(p)=1$
and $h^0(2p)=2$. Therefore, the natural map 
$H^0(\OO(2p))\to H^0(\OO(D+2p))$ is an isomorphism. Now
the fact that $r_{D+3p,p}$ is an isomorphism follows from the long exact sequence of cohomology
associated with the exact sequence
$$0\to\OO(D+2p)\to\OO(D+3p)\to \OO(3p)/\OO(2p)\to 0.$$
We have a natural direct sum decomposition 
$$H^0(\OO(D+3p)/\OO(p))\simeq H^0(\OO(D)/\OO)\oplus H^0(\OO(3p)/\OO(p))$$
and a boundary map
$$\de:H^0(\OO(D+3p)/\OO(p))\to H^1(\OO(p))\simeq H^1(\OO)/\lan\xi_p\ran.$$
We observe that the restriction of $\de$ to the summand $H^0(\OO(D)/\OO)$ is exactly
the isomorphism $\ov{\de_{D+p}}$, and the restriction of $\de$
to $H^0(\OO(3p)/\OO(p))$ is compatible with $\phi$.
Now start with a section $s\in H^0(\OO(D+3p))$ and write
$$s \mod\OO(p)=x+y$$ with
$x\in H^0(\OO(D)/\OO)$ and $y\in H^0(\OO(3p)/\OO(p))$. Then we have
$$\phi(r_{D+3p,p}(s))=\de(y),$$
$$\ov{\de_{D+p}}(r_{D+3p,D}(s))=\de(x).$$
Since $\de(x)+\de(y)=\de(s)=0$, our assertion follows.
\ed

\begin{cor}\label{4-Mas-prod-cor} 
Let $C$ be an irreducible projective curve with at most nodal singularities of arithmetic
genus $g\ge 2$,
and let $p\in C$ be a smooth point such that $H^0(C,\OO(2p))\not\sub H^0(C,\OO(p))$.
Then the Massey product \eqref{4-Mas-prod-eq} does not vanish.
\end{cor}

\Pf . By Proposition \ref{4-Mas-prod-prop}, we have to check that
$H^0(C,\OO(3p))=H^0(C,\OO(2p))$.
If $C$ is smooth then the divisor $2p$ is in the hyperelliptic system, and the assertion
follows easily. Thus, we can assume that $C$ is singular.
Let $\wt{C}\to C$ be the normalization of $C$,
so that $C$ is obtained by gluing pairs of distinct points $(a_i,b_i)$, $i=1,\ldots,s$, on $\wt{C}$.
We denote by $p\in\wt{C}$ the point corresponding to $p\in C$. If the genus of $\wt{C}$ is $\ge 2$
then it is hyperelliptic and the assertion follows as in the smooth case.

Now assume that $\wt{C}$ has genus $1$. The condition $h^0(C,\OO(2p))=2$ implies
that a nonconstant section $f\in H^0(\wt{C},\OO(2p))$ satisfies $f(a_i)=f(b_i)$ for $i=1,\ldots,s$.
Pick an element $h\in H^0(\wt{C},\OO(3p))\setminus H^0(\wt{C},\OO(2p))$.
Since the sections $(1,f,h)$ form a basis of $H^0(\wt{C},\OO(3p))$, they distinguish points of
$\wt{C}$, so we have $h(a_i)\neq h(b_i)$, and 
$h$ cannot descend to an element of $H^0(C,\OO(3p))$.
Hence, $H^0(C,\OO(3p))=H^0(C,\OO(2p))$ in this case. 

Finally, consider the case $\wt{C}=\P^1$. We can assume that $p=\infty$ and think of sections
of $\OO(np)$ on $\P^1$ as polynomials of degree $n$. Without loss of generality we can assume
that $b_i=-a_i$ for all $i$ (so that $t^2\in H^0(\P^1,\OO(2p))$ descends to a non-constant section of 
$\OO(2p)$ on $C$). Assume that there is a polynomial $h$ of degree $3$ such that
$h(a_i)=h(-a_i)$ for all $i$. Write $h=h_+ +h_-$, where $h_+$ is even and $h_-$ is odd.
Then we have $h_-(a_i)=0$ for $i=1,\ldots,s$. Since $h_-$ is an odd cubic polynomial, this implies
that $s=1$, which contradicts the assumption $g\ge 2$.
\ed

\begin{rem} There are examples of $(C,p)$, such that 
$H^0(C,\OO(3p))\not\sub H^0(C,\OO(2p))\not\sub H^0(C,\OO(p))$ and the arithmetic genus of
$C$ is $\ge 2$ (necessarily with $C$ reducible). The simplest example of genus $2$
is the union of two elliptic curves intersecting at one point ($p$ can be any smooth point).
\end{rem}

In the case when $C$ is a hyperelliptic smooth projective curve we can calculate
the Massey product \eqref{4-Mas-prod-eq} in terms of the corresponding ramification
points on $\P^1$.
Let $f:C\to\P^1$ be the morphism given by the hyperelliptic linear system,
so that $\OO(2p)\simeq f^*\OO(1)$. Let $p_1,\ldots,p_g$ be distinct Weierstrass
points on $C$ (i.e., ramification points of $f$), and let $a_i=f(p_i)\in\P^1$.
Then setting $D_i=\sum_{j\neq i}p_i$ we can use an isomorphism
$H^0(\OO(D_i)/\OO)\simeq H^1(\OO)/\lan X_i\ran$ (where $X_i=\xi_{p_i}$) 
and view the Massey product \eqref{4-Mas-prod-eq} for $p=p_i$ as
a map 
$$T_{p_i}^{\ot 3}\to H^0(\OO(D_i)/\OO)\simeq\bigoplus_{j\neq i} T_{p_j}.$$
Let $\a^{he}_{ij}: T_{p_i}^{\ot 3}\to T_{p_j}$ be the components of this map, where $i\neq j$.

\begin{prop} Assume $\cha(\k)\neq 2$.
Let $f: C\to\P^1$ be a double covering (with $C$ smooth) and let $p_1,\ldots,p_g$
be distinct ramification points with $f(p_i)=a_i\in\P^1$.
Let $\LL=\det^{-1}(f_*\OO_C)$ (which is a line bundle on $\P^1$ isomorphic to $\OO(g+1))$,
and let $s\in H^0(\P^1,\LL^2)$ be the canonical section vanishing at the ramification points of $f$.
Then we have canonical isomorphisms 
\begin{equation}\label{tangent-p-a-eq}
T_{p_i}\simeq\LL^{-1}|_{a_i}\ot T_{a_i}, 
\end{equation} 
where
$T_{a_i}$ is the tangent space to $\P^1$ at $a_i$. Let us write
$$\a^{he}_{ij}=\b_{ij}\cdot \kappa_i,$$
where $\kappa_i\in \LL^2|_{a_i}\ot T_{p_i}^{-1}$ is the natural trivialization given by $s$,
and 
$$\b_{ij}\in \LL|_{a_i}\ot T_{a_i}^{-2}\ot \LL_{a_j}^{-1}\ot T_{a_j}.$$
Assume that $a_i\in \A^1=\P^1\setminus\{\infty\}$ for $i=1,\ldots,g$. Then using
a trivialization of $\LL|_{\A^1}$ and the natural trivialization of $T_{\P^1}|_{\A^1}$, we have
$$\b_{ij}=\frac{1}{a_j-a_i}\cdot\prod_{k\neq i,j}\frac{a_i-a_k}{a_j-a_k}.$$
Hence, using these trivializations, and in addition an isomorphism $\LL^2\simeq \OO_{\P^1}(2g+2)$, 
we get
$$\a^{he}_{ij}=\frac{F'(a_i)}{a_j-a_i}\cdot\prod_{k\neq i,j}\frac{a_i-a_k}{a_j-a_k},$$
where $F$ is the polynomial corresponding to the section $s$.
\end{prop}

\Pf . We can realize $C$ as the relative spectrum of the sheaf of $\OO_{\P^1}$-algebras $\AA=f_*\OO_C$.
Since $\cha(\k)\neq 2$, we have a canonical isomorphism
$$\AA=\OO_{\P^1}\oplus \LL^{-1},$$
where $\LL\simeq\OO_{\P^1}(g+1)$ and the product $\LL^{-1}\ot\LL^{-1}\to \OO$ is given by a
section $s\in H^0(\P^1,\LL^2)$ that has simple zeros at the $2g+2$ ramification points of $f$.
Recall that by Proposition \ref{4-Mas-prod-prop}(ii), to compute $\a^{he}_{ij}$
we have to choose an element $h_i\in H^0(\OO_C(3p_i+\sum_{j\neq i}p_j))$, such that
it projects to a given generator of $T_{p_i}^{\ot 3}\simeq \OO_C(3p_i)/\OO_C(2p_i)$, and
then consider the polar part of $h_i$ at $p_j$.
Note that the maximal ideal sheaf in $\AA$ corresponding to $p_i\in C$ is 
$$\mg_{p_i}=\mg_{a_i}\oplus \LL^{-1}\sub \OO_{\P^1}\oplus\LL^{-1},$$
where $\mg_a$ is the maximal ideal sheaf corresponding to $a_i=f(a_i)\in\P^1$.
Hence, 
$$f_*\OO_C(p_i)=\mg^{-1}_{p_i}\simeq \OO_{\P^1}\oplus \mg^{-1}_{a_i}\LL^{-1} \text{ and }$$
$$T_{p_i}\simeq \mg^{-1}_{p_i}/\OO_C\simeq \mg^{-1}_{a_i}\LL^{-1}/\LL^{-1}\simeq\LL^{-1}|_{a_i}\ot T_{a_i}.$$
We define \eqref{tangent-p-a-eq} to be this latter isomorphism.
On the other hand, we have a natural isomorphism
$$f_*\OO_C(2p_i)\simeq \AA(a_i)=\OO_{\P^1}(a_i)\oplus \LL^{-1}(a_i)$$
such that the induced isomorphism
$$T_{p_i}^2\simeq \OO_C(2p)/\OO_C(p)\simeq \mg^{-1}_{a_i}/\OO_{\P^1}\simeq T_{a_i}$$
differs from the square of \eqref{tangent-p-a-eq} by $\kappa_i\in\LL^2|_{a_i}\ot T_{p_i}^{-1}$.
Similarly, we get a natural isomorphism
\begin{equation}\label{hyper-rat-fun-eq}
f_*\OO_C(3p_i+\sum_{j\neq i}p_j)\simeq \OO_{\P^1}(a_i)\oplus \LL^{-1}(2a_i+\sum_{j\neq i}a_j),
\end{equation}
compatible with isomorphisms
$$T_{p_j}\simeq\OO_C(p_j)/\OO_C\simeq \LL^{-1}(a_j)|_{a_j}\simeq \LL^{-1}|_{a_j}\ot T_{a_j},$$
$$\OO_C(3p_i)/\OO_C(2p_i)\simeq \LL^{-1}(2a_i)|_{a_i}\simeq \LL^{-1}|_{a_i}\ot T_{a_i}^{\ot 2},$$
of which the first is \eqref{tangent-p-a-eq} and the second differs from \eqref{tangent-p-a-eq} by $\kappa_i$.
In terms of the isomorphism \eqref{hyper-rat-fun-eq} we can find the required rational function by choosing
a global section $\wt{h}_i$ of $\LL^{-1}(2a_i+\sum_{j\neq i}a_j)$. Then, if $\wt{h}_i$ projects to a given generator
of $\LL^{-1}(2a_i)|_{a_i}\simeq \LL^{-1}|_{a_i}\ot T_{a_i}^2$, the polar part of $\wt{h}_i$ at $a_j$
(which is an element of $\LL^{-1}|_{a_j}\ot T_{a_j}$) will give $\b_{ij}$.

Now assume that all $a_j\in\A^1=\P^1\setminus\{\infty\}$ and let $t$ be the natural coordinate on $\A^1$.
Let us fix a section $u\in H^0(\A^1,\LL^{-1})$ that gives a trivialization of $\LL^{-1}$ on $\A^1$ with the pole of order $g+1$ at infinity. Then $\frac{u}{(t-a_i)^2}$ (resp., $\frac{u}{t-a_j}$) induces a trivialization of 
$\LL^{-1}|_{a_i}\ot T_i^2$ (resp., $\LL^{-1}|_{a_j}\ot T_j$).
Thus, we can set
$$\wt{h}_i=\frac{u}{(t-a_i)^2}\cdot\prod_{j\neq i}\frac{a_i-a_j}{t-a_j}.$$
Taking the residue of $\wt{h}_i$ at $a_j$ we get the desired formula for $\b_{ij}$.
\ed

\subsection{Consequences for the $A_{\infty}$ structure}

Let $(C,p_1,\ldots,p_n)$ be a smooth projective curve of genus $g\ge 1$ with $n$ marked points (where $n\ge 1$), such  that $h^0(p_1+\ldots+p_n)=1$.
Let $E=E_{g,n}$ be the $\Ext$-algebra of the generator
$\OO_C\oplus\OO_{p_1}\oplus\ldots\oplus\OO_{p_n}$ of $D^b(C)$.
By the homological perturbation theory, we have a minimal $A_{\infty}$-structure
on $E$ extending the associative product on $E$, defined uniquely up to
$A_\infty$-equivalence.

\begin{thm}\label{m3-m4-thm} 
The $A_{\infty}$-structure on $E$ coming from the data $(C,p_1,\ldots,p_n)$ is equivalent
to the one with $m_3=0$ if and only if either $g=1$ or
$C$ is hyperelliptic and $p_1,\ldots,p_n$
are Weierstrass points. If $m_3=0$ and $g>1$ then $m_4$ is always nontrivial.
\end{thm}

\Pf . Assume first that $n=g$.
A minimal $A_{\infty}$-structure is equivalent to the one with $m_3=0$ if and only
if the Hochschild cohomology class given by $m_3$ is trivial.
By Proposition \ref{Hoch-prop}, this happens exactly when
$$m_3(B_i,Y_i,A_i)\in\lan X_i\ran$$
for all $i=1,\ldots,n$.
By Corollary \ref{mn-Massey-cor}, this is equivalent to the vanishing of the Massey products
\eqref{Mas-prod-eq} for $p=p_1,\ldots,p_n$.
Now Proposition \ref{van-Mas-prop} tells that this is equivalent to $C$ being hyperelliptic
and $p_1,\ldots,p_n$ being Weierstrass points.

In the case $n<g$ considering the same Massey products shows that the condition for
$C$ to be hyperelliptic and for $p_1,\ldots,p_n$ to be Weierstrass points is necessary. 
Conversely, if we have such $n$-tuple of Weierstrass points on a hyperelliptic curve we can 
complete it to a $g$-tuple of Weierstrass points $p_1,\ldots, p_g$ still satisfying
the condition $h^0(p_1+\ldots+p_g)=1$ (see Lemma \ref{W-points-lem} below).
By the first part of the proof, the $A_{\infty}$-structure on $\OO, \OO_{p_1},\ldots,\OO_{p_g}$
can be chosen to have trivial $m_3$, as required.

The second assertion follows from the nontriviality of the quadruple Massey product
\eqref{4-Mas-prod-eq} for a Weierstrass point on a hyperelliptic curve 
(see Corollary \ref{4-Mas-prod-cor}). To connect this Massey product to $m_4$ we use
Corollary \ref{mn-Massey-cor}, noting that the needed dg Massey product is defined by
Lemma \ref{4Mas-def-lem}.
\ed

\begin{lem}\label{W-points-lem} Let $p_1,\ldots,p_n$ be distinct Weierstrass points
on a hyperelliptic curve, where $n\le g$. Then $h^0(p_1+\ldots+p_n)=1$.
\end{lem}

\Pf . It is enough to consider the case $n=g$, in which case we have to check that
$h^1(D)=h^0(K-D)=0$, where $D=p_1+\ldots+p_g$. Indeed, otherwise we would have
$K=D+D'$ for some effective divisor $D'$ of degree $g-2$. Since every effective canonical
divisor on $C$ is a sum of $g-1$ fibers of the hyperelliptic map $f:C\to\P^1$, this would
imply that $f(D+D')$ is supported at $\le g-1$ points, which is a contradiction.
\ed

\section{Rational functions on $\MM_{g,g}$ associated with Massey products}\label{rat-sec}

\subsection{Triple products as sections of line bundles over the moduli spaces}
\label{sections-moduli-sec}

Let $C$ be a projective curve of arithmetic genus $g\ge 2$, and let 
$p_1,\ldots, p_g$ be distinct smooth points such that $h^0(p_1+\ldots+p_g)=1$.
Then, by Proposition \ref{Hoch-prop} and Corollary \ref{alpha-for-cor},
the Hochschild class of $m_3$ on $E_{g,g}$ (where the $A_{\infty}$-structure comes from 
$(C,p_1,\ldots,p_g)$) is determined
by the collection of elements $\a_{ij}\in\Hom_k(T_{p_i}^{\ot 2}, T_{p_j})$, $i\neq j$, given by
\begin{equation}\label{alpha-def-eq}
\a_{ij}=r_{D_i+2p_i,p_j}\circ r_{D_i+2p_i,p_i}^{-1},
\end{equation} 
where we use the restriction maps
\begin{eqnarray}\label{restriction-maps-eq}
r_{D_i+2p_i,p_i}:H^0(\OO(D_i+2p_i))/H^0(\OO)\rTo{\sim} H^0(\OO(2p_i)/\OO(p_i))\simeq T_{p_i}^{\ot 2}
\nonumber \\
r_{D_i+2p_i,p_j}:H^0(\OO(D_i+2p_i))/H^0(\OO)\to H^0(\OO(p_j)/\OO)\simeq T_{p_j},
\end{eqnarray}
where $D_i=\sum_{j\neq i} p_j$.
In particular, this construction makes sense over the open substack $\UU\sub\ov{\MM}_{g,g}$ 
of the Deligne-Mumford stack of stable curves with marked points, 
corresponding to $(C,p_1,\ldots,p_g)$ such that
$h^0(p_1+\ldots+p_g)=1$. Thus, $\a_{ij}$ can be viewed as a section over $\UU$ of
the line bundle $L_i^2\ot L_j^{-1}$, where $L_i:=p_i^*K$ (the pullback of the relative canonical
class on the universal curve).

Let us set $D_{ij}=\sum_{m\neq i,j} p_m$.
The zero locus of $\a_{ij}$ is supported on the divisor $Z_{ij}\sub\ov{\MM}_{g,g}$ of 
$(C,p_1,\ldots,p_g)$ such that $h^0(2p_i+D_{ij})>1$.
In particular, $\a_{ij}$ is nonzero.
The complement to $\UU$
is the divisor $Z\sub\ov{\MM}_{g,g}$ of $(C,p_1,\ldots,p_g)$ such that 
$h^0(p_1+\ldots+p_g)>1$. More precisely, we define $Z$ as the degeneration locus
of the map $H^0(\OO(p_1+\ldots+p_g)/\OO)\to H^1(\OO)$, which
is the zero locus of a section of the line bundle $\det(\La)^{-1}\ot L_1\ot\ldots\ot L_g$
on $\ov{\MM}_{g,g}$, where $\La$ is the Hodge bundle.
Similarly, $Z_{ij}$ is defined as the degeneration locus of the map
$H^0(\OO(2p_i+D_{ij})/\OO)\to H^1(\OO)$, so it is the zero locus of a section of
$\det(\La)^{-1}\ot L_i^3\ot\bigotimes_{m\neq i,j} L_m$.  

Note that the divisors $Z$ and $Z_{ij}$ have in general many irreducible components
(they contain some boundary components). 

\begin{prop}\label{moduli-section-prop}
The section $\a_{ij}\in\Ga(\UU,L_i^2\ot L_j^{-1})$ extends to a global section 
$$\wt{\a}_{ij}\in \Ga(\ov{\MM}_{g,g},L_i^2\ot L_j^{-1}(Z))$$
such that the zero locus of $\wt{\a}_{ij}$ is exactly $Z_{ij}$. 
\end{prop}

The proof will be based on the following general fact from tensor algebra.
Recall that for a morphism of vector bundles $\phi:V\to W$ such that $r=\rk W=\rk V-1$
one has a canonical map
$$k_\phi:\det(V)\ot \det(W^*)\to V$$
such that $\phi\circ k_\phi=0$. Namely, $k_\phi$ is obtained by tensoring with
$\det(V)$ from the map
$${\bigwedge}^r(\phi^*):\det(W^*)\to {\bigwedge}^r(V^*)$$
using the natural isomorphism $\det(V)\ot\bigwedge^r(V^*)\simeq V$.

\begin{lem}\label{det-lem}
Let $0\to V_1\rTo{\iota} V\rTo{\pi} L\to 0$ be an exact sequence of vector bundles,
where $L$ is a line bundle, and let $\phi:V\to W$ be a morphism of vector
bundles, where $r=\rk W=\rk V_1=\rk V-1$. 
Let $Z\sub S$ be the degeneration divisor of the restriction
$\phi_1=\phi|_{V_1}:V_1\to W$. Then $Z$ coincides with the
vanishing locus of the composed map
$$\det(V)\ot\det(W^*)\rTo{k_\phi} V \rTo{\pi} L.$$
\end{lem}

\Pf . Note that $Z$ is the vanishing locus of $\det(\phi_1):\det(V_1)\to\det(W)$,
or equivalently, of the dual map $\det(\phi_1^*)$.
Thus, the assertion follows from the commutativity of the diagram
\begin{diagram}
\det(V)\ot\det(W^*)&\rTo{\bigwedge^r(\phi^*)}&\det(V)\ot{\bigwedge}^r(V^*)&\rTo{\sim}& V\\
&\rdTo{\det(\phi_1^*)}&\dTo{\bigwedge^r(\iota^*)}&&\dTo{\pm\pi}\\
&&\det(V)\ot\det(V_1^*)&\rTo{\sim}& L
\end{diagram}
since the composition of arrows in the top row is $k_\phi$.
\ed

\noindent
{\it Proof of Proposition \ref{moduli-section-prop}.}
Let $V$ be the bundle on $\ov{\MM}_{g,g}$ with the fiber 
$H^0(C,\OO(2p_i+D_i)/\OO)$ over $(C,p_1,\ldots,p_g)$, and let $W=\La^*$,
so the fiber of $W$ at $(C,p_1,\ldots,p_g)$ is $H^1(C,\OO)$. We have
a natural connecting homomorphism $\phi:V\to W$. We have natural restriction maps
$\pi_i:V\to L_i^{-2}$ and $\pi_j:V\to L_j^{-1}$. Applying Lemma \ref{det-lem}
to the exact sequence
of bundles 
$$0\to V'\to V\to L_i^{-2}\to 0,$$ where
$V'$ is the bundle on $\ov{\MM}_{g,g}$ with the fiber $H^0(C,\OO(p_i+D_i))$,
we see that the divisor $Z\sub\ov{\MM}_{g,g}$ 
coincides with the vanishing locus of the composition
$$\det(V)\ot\det(W^*)\rTo{k_\phi} V \rTo{\pi_i} L_i^{-2}.$$
Note that over $\UU$ the image of $k_\phi$ generates $\ker(\phi)$, and
$\ker(\phi)$ is a bundle with the fiber
$$\ker(H^0(C,\OO(2p_i+D_i)/\OO)\to H^1(C,\OO))\simeq H^0(C,\OO(2p_i+D_i))/H^0(C,\OO).$$
Thus, we can replace the restriction maps $r_{D_i+2p_i,p_i}$ and $r_{D_i+2p_i,p_j}$ used in defining $\a_{ij}$
(see \eqref{restriction-maps-eq}) with the morphisms
$\pi_i\circ k_\phi$ and $\pi_j\circ k_\phi$, respectively.
Since $\pi_i\circ k_\phi$ induces an isomorphism
$$\det(V)\ot\det(W^*)\simeq L_i^{-2}(-Z),$$
we obtain the global morphism
$$L_i^{-2}(-Z)\simeq\det(V)\ot\det(W^*)\rTo{\pi_j\circ k_\phi} L_j^{-1}$$
which gives the required global section $\wt{\a}_{ij}\in \Ga(\ov{\MM}_{g,g}, L_i^2\ot L_j^{-1}(Z))$.
Now applying Lemma \ref{det-lem} to the exact sequence
$$0\to V''\to V\to L_j^{-1}\to 0,$$
where $V''$ is the bundle on $\ov{\MM}_{g,g}$ with the fiber $H^0(C,\OO(2p_i+D_{ij}))$,
we see that the vanishing locus of $\pi_j\circ k_\phi$ is exactly $Z_{ij}$.
\ed

\begin{rem} It is not essential to work with stable curves in the above argument. The result
similar to Proposition \ref{moduli-section-prop} would work with other modular compactifications
of $\MM_{g,g}$.
\end{rem}

\subsection{Rational functions}\label{rat-fun-sec}

Let $\MM^{(1)}_{g,g}\to \MM_{g,g}$ (resp.,
$\ov{\MM}^{(1)}_{g,g}\to\ov{\MM}_{g,g}$) be the $\G_m^g$-torsor 
corresponding to choices of nonzero tangent vectors at each of the marked points.
Then the line bundles $L_i$ are naturally trivialized on $\MM^{(1)}_{g,g}$, so
we can view each section $\a_{ij}$ as a rational function on $\MM^{(1)}_{g,g}$.
This gives a rational map
\begin{equation}\label{big-rat-map-eq}
\a:\MM^{(1)}_{g,g}\rTo{(\a_{ij})}\G_m^{g^2-g}.
\end{equation}
On the other hand,
considering rational monomials in $\a_{ij}$, we can get rational functions on $\MM_{g,g}$.
Namely, consider the homomorphism of groups 
$$\varphi:\Z^{g^2-g}\to\Z^g: e_{ij}\to 2e_i-e_j,$$
where $\Z^{g^2-g}$ (resp., $\Z^g$) has a basis $(e_{ij})_{i\neq j}$ (resp., $e_i$), where $i,j\le g$.
Then for every element $x=\sum n_{ij}e_{ij}\in\ker(\varphi)$ the expression
$$\a^x:=\prod \a_{ij}^{n_{ij}}$$
is a rational function on $\MM_{g,g}$.
It is easy to see that $\ker(\varphi)$ has rank $g^2-2g$, so choosing a basis
$b_1,\ldots,b_{g^2-2g}$ in $\ker(\varphi)$ we obtain a rational map
\begin{equation}\label{rat-map-eq}
\ov{\a}:\MM_{g,g}\rTo{{\a^{b_1}},\ldots,\a^{b_{g^2-2g}}}\G_m^{g^2-2g}.
\end{equation}
Note that the rational map \eqref{big-rat-map-eq} is $\G_m^g$-equivariant, where $(\la_1,\ldots,\la_g)$
acts on $\G_m^{g^2-g}$ via the homomorphism $\varphi^*:\G_m^g\to\G_m^{g^2-g}$,
dual to $\varphi$, and the map $\ov{\a}$ can be viewed as the induced rational map
of quotients by $\G_m^g$.

\begin{thm}\label{birational-thm} Let $\cha(\k)=0$.
If $g\ge 6$ then the map \eqref{rat-map-eq} is birational
onto its image. 
\end{thm}

The proof of this theorem will be given in Section \ref{recon-sec}. 
The result is optimal, since for $g\le 5$ we have $\dim\MM_{g,g}>g^2-2g$.
In fact, for $g\le 5$ the map \eqref{rat-map-eq} is dominant (see Theorem \ref{dominant-thm}
below).


\begin{prop}\label{image-prop} 
Let $g\ge 3$. For a generic curve $C$ the restriction of $\ov{\a}$
gives a rational map 
$$\ov{\a}_C: C^g\to\G_m^{g^2-2g}$$ 
with generically injective tangent map.
Hence, the image of $\ov{\a}_C$ has dimension $g$.
\end{prop}

\Pf . Using the sections $\a_{ij}$ on $\UU\sub\ov{\MM}_{g,g}$ (see Section
\ref{sections-moduli-sec}) we can extend the map $\ov{\a}$ to stable curves. 
It is enough to construct 
a stable curve $(C,p_1,\ldots,p_g)$ in $\UU$ for which the assertion is true.
Let us consider the wheel of $\P^1$'s with 
$g$ components $C_1,\ldots,C_g$, so that $1\in C_i$ is glued to $0\in C_{i+1}$
(we think of indices as elements of $\Z/g\Z$). Now consider the nodal 
curve $C$ obtained as the union of this wheel with 
one more component $C_{\infty}\simeq \P^1$ which intersects each component $C_i$ at one
point $\infty\in C_i$ (we fix the corresponding $g$ distinct points on $C_{\infty}$).
Note that the arithmetic genus of $C$ is $g$.
We choose marked points $p_1,\ldots,p_g$, so that $p_i=\la_i\in C_i\setminus\{0,1,\infty\}$.

Let us compute $\a_{1i}$. By definition, for this we have to produce a non-constant element 
$f\in H^0(C,\OO(2p_1+p_2+\ldots+p_g))$. Such a function is given by a collection of functions
$(f_1,\ldots,f_g,f_{\infty})$, where $f_1\in H^0(C_1,\OO(2p_1))$,
$f_i\in H^0(C_i,\OO(p_i))$ for $i=2,\ldots,g$ and $f_{\infty}$ is a constant, subject to
the constraints
$$f_i(1)=f_{i+1}(0), $$
$$f_i(\infty)=f_{\infty},$$
where $i=1,\ldots,g$.
Subtracting a constant from $f$ we can assume that $f_{\infty}=0$.
Then we can take 
$$f_1(t)=\frac{1}{(t-\la_1)^2}+\frac{y_1}{t-\la_1},$$
$$f_i(t)=\frac{y_i}{t-\la_i},$$
for some constants $y_1,\ldots,y_g$,
and the equations become
$$\frac{1}{(1-\la_1)^2}+\frac{y_1}{1-\la_1}=-\frac{y_2}{\la_2},$$
$$\frac{y_i}{1-\la_i}=-\frac{y_{i+1}}{\la_{i+1}}, \ i=2,\ldots,g-1,$$
$$\frac{y_g}{1-\la_g}=\frac{1}{\la_1^2}-\frac{y_1}{\la_1}.$$
Solving this system we obtain 
$$\a_{12}=y_2=\frac{\la_2}{\la_1(\la_1-1)^2(a-1)}, \ 
\a_{13}=y_3=\frac{\la_2\la_3}{\la_1(\la_1-1)^2(\la_2-1)(a-1)}, \ \text{ etc.},$$
where 
$$a=\frac{\la_1\la_2\ldots\la_g}{(\la_1-1)(\la_2-1)\ldots(\la_g-1)}.$$

Now we find
$$\frac{\a_{i,i+1}^2\a_{i+1,i+3}}{\a_{i,i+2}^2\a_{i+2,i+3}}=\frac{\la_{i+2}-1}{\la_{i+1}} \
\text{ for } i=1,\ldots,g,$$
which implies that the parameters $\la_1,\ldots,\la_g$ can be recovered from
the image of the map \eqref{rat-map-eq}.
\ed


\begin{ex} In the case $g=2$ the homomorphism $\varphi^*:\G_m^2\to \G_m^{2}$ has 
kernel $\Z/3\Z\sub \G_m^2$, generated by $(\zeta_3,\zeta_3^{-1})$, where
$\zeta_3$ is a primitive $3$rd root of unity. Hence, the map $\a$ in this case factors
through a rational map
$$\a':\MM^{(1)}_{2,2}/(\Z/3\Z)\to\G_m^2,$$
and the induced map
\begin{equation}\label{genus-2-birat-map}
\MM^{(1)}_{2,2}/(\Z/3\Z)\rTo{(\pi,\a')} \MM_{2,2}\times\G_m^2
\end{equation}
is birational (where $\pi$ is the natural projection to $\MM_{2,2}$).
More explicitly,
the inverse rational map to \eqref{genus-2-birat-map} sends $(C,p_1,p_2,\la,\mu)$ to
$(C,p_1,p_2,v_1,v_2)$, where the tangent vectors  $v_1\in T_{p_1}$, $v_2\in T_{p_2}$
are defined uniquely
up to the $\Z/3\Z$-action (generated by $(v_1,v_2)\mapsto (\zeta_3 v_1,\zeta_3^{-1} v_2)$)
 by the condition that there exist rational functions $f_1\in H^0(C,\OO(2p_1+p_2))$,
$f_2\in H^0(C,\OO(p_1+2p_2))$ with
\begin{align*}
& f_1\equiv v_1^2\mod\OO(p_1+p_2), \ f_1\equiv \la v_2\mod\OO(2p_1),\\
& f_2\equiv v_2^2\mod\OO(p_1+p_2), \ f_2\equiv \mu v_1\mod\OO(2p_2).
\end{align*}
\end{ex}

\begin{ex} In the case $g=3$ the space $\MM^{(1)}_{3,3}$ is $12$-dimensional. By
Proposition \ref{image-prop}, for a generic curve $C$ of genus $3$ the rational map
$$\ov{\a}_C: C^3\to \G_m^3$$
is generically \'etale. Hence, 
at generic point of $\MM^{(1)}_{3,3}$ the fibers of
the two dominant (rational) maps to $6$-dimensional spaces
\begin{diagram}
\MM^{(1)}_{3,3}&\rTo{\a}&\G_m^6\\
\dTo{\pi}\\
\MM_3
\end{diagram}
are transversal.
\end{ex}

\subsection{Interpretation in terms of tangent lines}\label{geom-int-sec}

Let $C$ be a smooth projective curve of genus $g\ge 2$.
Let $L$ be a base point free line bundle on $C$. For a point
$p\in C$ let $\ev_p\in L|_p\ot H^0(C,L)^*$ denote the functional of evaluation at $p$.
Then the tangent map at a point $p\in C$ to the
map 
$$\varphi_L:C\rTo{|L|}\P(H^0(C,L)^*),$$ 
given by the linear system $|L|$, is the map 
$$T_pC\to L|_p\ot H^0(C,L)^*/\lan\ev_p\ran\simeq L|_p\ot H^0(C,L(-p))^*,$$
dual to the evaluation functional for $L(-p)$,
$$H^0(C,L(-p))\to L(-p)|_p\simeq (T_pC)^*\ot L|_p.$$
In the case of the canonical line bundle $L=\om_C$, under the duality
$H^0(C,\om_C)^*\simeq H^1(C,\OO_C)$ the functional $\ev_p$ corresponds to
the element $\xi_p\in H^1(C,\OO_C)$, obtained from the connecting homomorphism
$H^0(\OO(p)/\OO)\to H^1(\OO)$. Hence, the tangent map to the canonical morphism
$\varphi_{\om_C}:C\to\P(H^0(C,\om_C)^*)$ at $p\in C$ can be identified with
the connecting homomorphism
$$\de'_p:
T_pC\simeq T_p^*C\ot H^0(C,\OO(2p)/\OO(p))\to  T_p^*C\ot H^1(C,\OO(p))\simeq 
T_p^*C\ot H^1(C,\OO)/\lan\xi_p\ran,$$ 
which is exactly the triple Massey product considered in Section \ref{3Mas-sec}.

Now recall that for $g$ distinct points $p_1,\ldots,p_g\in C$ such that $h^0(p_1+\ldots+p_g)=1$
the maps $\a_{ij}:T_{p_i}^2\to T_{p_j}$ for $i\neq j$, considered above, 
can be identified with the components of the same Massey product for $p=p_i$, 
up to a sign (see Proposition
\ref{3-Mas-prod-prop} and Corollary \ref{alpha-for-cor}).
This leads to the following identification of the rows of the matrix $(-\a_{ij})$ with
the coordinates of the tangent map to the morphism given by $|\om_C|$.

\begin{prop}\label{canonical-tangent-line-prop}
The components of the tangent map $\de'_{p_i}$ to $\varphi_{\om_C}$ at $p_i$, 
with respect to the decomposition $H^1(C,\OO)\simeq\bigoplus_j T_{p_j}$,
are given by tensoring with $-\a_{ij}$.
\end{prop}

Note that the position of the tangent line to $C$ at $p_i$ in $\P^{g-1}$ is
recorded by $(\a_{ij})$ with fixed $i$, viewed as homogeneous coordinates.
In order to recover the same data as the map $\ov{\a}$, we note that there is a canonical
identification of the tangent line to $C$ at $p_i$ with the fiber of the tautological line
bundle $\OO_{\P^{g-1}}(-1)$ at $p_i$. Thus, Theorem \ref{birational-thm} leads to the following result.

\begin{cor}\label{canonical-tangents-cor}
Let $\cha(\k)=0$ and $g\ge 6$.
Let us associate with generic $(C,p_\bullet)\in\MM_{g,g}$ the collection
(for $i=1,\ldots,g$) of points 
$x_i=\varphi_{\om_C}(p_i)\in\P^{g-1}$
and of tangent lines $L_i\sub\P^{g-1}$ to $\varphi_{\om_C}(C)$ at $x_i$
together with identification of each tangent space
$T_{x_i}L_i$ with the fiber of the tautological line bundle $\OO(-1)|_{x_i}$.
Then generic $(C,p_\bullet)$ can be recovered from these data (viewed up to projective
transformations).
\end{cor}

Next, consider the map $C\to\P^g$ given by the linear system $|2D|$, where $D=\sum_{i=1}^g p_i$.
If $h^0(2D-K)=0$ (which is true generically) then this map is an embedding and its image
is a degree $2g$ curve in $\P^g$.
Note that the section $1\in H^0(C,\OO(2D))$ corresponds to a hyperplane $H\sub\P^g$
which is tangent to $C$ at all $g$ points $p_1,\ldots,p_g$.
Also the condition $h^0(D)=1$ means that $p_1,\ldots,p_g$ are in general position in $H$.

Now suppose we are given any degree $2g$ curve 
$C\sub\P^g$, and a linear form $\ell\in H^0(\P^g,\OO(1))$
such that the corresponding hyperplane $H=(\ell=0)$ is tangent to $C$ at $g$ points
$p_1,\ldots,p_g$ that are smooth points of $C$ and are in general linear position (we
assume also that $C\not\sub H$).
Let $L=\OO(1)|_C$. Since $\deg(L)=2g$ and the section $\ell$ vanishes along the divisor
$2p_1+\ldots+2p_g$, it induces an isomorphism 
\begin{equation}\label{tautological-isom-eq}
\OO(1)|_{p_i}\simeq \OO_C(2p_i)/\OO(p_i)\simeq T_{p_i}^2
\end{equation}
for each $i$. 
Since $p_1,\ldots,p_g$ are in general position, we obtain an isomorphism
$$H^0(H,\OO(1))\simeq\bigoplus_{i=1}^g\OO(1)|_{p_i}\simeq\bigoplus_{i=1}^g T_{p_i}^2.$$
Therefore, we have a canonical isomorphism
\begin{equation}\label{tangent-space-decomposition}
T_{p_i}H\simeq \bigoplus_{j\neq i} T_{p_i}^2\ot T_{p_j}^{-2}.
\end{equation}
Now, since $H$ is tangent to $C$ at each $p_i$,
the tangent map to the embedding $C\to\P^g$ at $p_i$ gives a linear map
\begin{equation}\label{geom-alpha-eq}
T_{p_i}\to T_{p_i}H\simeq \bigoplus_{j\neq i} T_{p_i}^2\ot T_{p_j}^{-2}.
\end{equation}
This time the tangent map will be given by a column of the matrix $(\a_{ij})$.

\begin{prop}\label{tangent-line-prop} Suppose $C\hra\P^g$ is a smooth projective curve embedded
by the linear system $|2D|$, where $D=\sum_{i=1}^g p_i$ and $h^0(D)=1$. Also, let
$\ell$ be the section $1\in H^0(C,\OO(2D))\simeq H^0(\P^g,\OO(1))$.
Then the components
of the map \eqref{geom-alpha-eq} are given by tensoring with $\a_{ji}$ (see \eqref{alpha-def-eq}).
\end{prop}

\Pf . We have $L=\OO(1)|_C=\OO_C(2D)$.
The point $p_i$ corresponds to the functional
$$H^0(\P^g,\OO(1))=H^0(C,L)\to L|_{p_i},$$
and the tangent map to the embedding $C\to \P^g$ at $p_i$
corresponds to the dual of the natural restriction map
\begin{equation}\label{tangent-curve-embedding-eq}
H^0(C,L(-p_i))/(1)\to L(-p_i)|_{p_i}\simeq L|_{p_i}\ot (T_{p_i}C)^*.
\end{equation}
The components of the direct sum decomposition \eqref{tangent-space-decomposition} 
are exactly the subspaces
$$H^0(C,\OO(2p_i+D_i))/(1)\sub H^0(C,L)/(1)\simeq H^0(H,\OO(1)).$$
The subspace $H^0(C,L(-p_i))/(1)\sub H^0(C,L)/(1)$ is the direct sum of the components
$H^0(C,\OO(2p_j+D_j))\simeq (T_{p_j}C)^{\ot 2}$ for $j\neq i$. Furthermore, the restriction
of \eqref{tangent-curve-embedding-eq} to the subspace
$H^0(C,\OO(2p_j+D_j))$ is exactly $\a_{ji}$. This immediately implies the assertion.
\ed

Similarly to Corollary \ref{canonical-tangents-cor} this leads to the following result.

\begin{cor}\label{another-tangents-cor} 
Let $\cha(\k)=0$ and $g\ge 6$. Then a generic  
$(C,p_1,\ldots,p_g)\in\MM_{g,g}$ is uniquely determined by the configuration of
$g$ points $p_1,\ldots,p_g$ and $g$ tangent lines $L_i$ to $C$ at these points 
in the embedding given by the linear system
$|2(p_1+\ldots+p_g)|$, together with identifications $(T_{p_i}L_i)^2\simeq\OO(1)|_{p_i}$
obtained from \eqref{tautological-isom-eq}.
\end{cor}

\begin{rem} One can ask whether in Corollaries \ref{canonical-tangents-cor} and
\ref{another-tangents-cor} it is enough to consider simply the configuration of points
$(p_i)$ and lines $(L_i)$, for sufficiently large $g$. We do not know the answer. Note
that the number of parameters describing such a configuration is
$g^2-3g$, so one should take $g\ge 7$ in order for this to have a chance to be true.
\end{rem}

Recall that for a line bundle $L$ one has the Wahl map (see \cite{Wahl})
$$W_L:{\bigwedge}^2 H^0(C,L)\to H^0(C,L^2\ot\om_C)$$
By definition, 
$$W_L(s_1\wedge s_2)(p)=\varphi'_1(p)\varphi_2(p)-\varphi'_2(p)\varphi_1(p),$$
where $\varphi_1,\varphi_2$ are local functions at $p$ corresponding to $s_1,s_2$ via some local trivialization of $L$. In invariant terms, the functional 
$$W_{L,p}:{\bigwedge}^2 H^0(C,L)\to (L^2\ot\om_C)|_p: s_1\wedge s_2\to 
W_L(s_1\wedge s_2)(p)$$ 
is given by restricting to a neighborhood $U\sub C$ of $p$ and applying the composition
$${\bigwedge}^2H^0(U,L)\to H^0(U,L(-p))\ot L|_p\rTo{\ev_p\ot\id} 
L(-p)|_p\ot L|_p\simeq (L^2\ot\om_C)|_p,$$
where the first map is induced by the exact sequence
$$0\to H^0(U,L(-p))\to H^0(U,L)\to L|_p\to 0.$$ 
In the case when $L$ is base point free we can take $U=C$, and we see that 
$W_{L,p}$ is essentially given by the Pl\"ucker coordinates of the tangent line to
$\varphi_L(C)\sub \P(H^0(C,L)^*)$ at $p$. 

Now let $(C,p_1,\ldots,p_g)$ be as before. 
Given the interpretation of $(\a_{ij})$ in terms of tangent lines from Propositions
\ref{canonical-tangent-line-prop} and \ref{tangent-line-prop}, we can relate it to
the Wahl maps $W_L$ associated with $L=\om_C$ and $L=\OO_C(2D)$.

In the case $L=\om_C$ we have a natural
decomposition $H^0(C,\om_C)=\bigoplus_{i=1}^g T_{p_i}^*$,
and for $i\neq j$ the restriction
$$T_{p_i}^{-1}\ot T_{p_j}^{-1}\hra {\bigwedge}^2 H^0(C,\om_C)\rTo{W_{\om_C,p_i}} T_{p_i}^{-3}$$
is given by tensoring with $-\a_{ij}$. This completely determines $W_{\om_C,p_i}$ since
its restrictions to $T_{p_j}^{-1}\ot T_{p_k}^{-1}\sub {\bigwedge}^2 H^0(C,\om_C)$ are zero
for $j\neq i$, $k\neq i$.

In the case $L=\OO_C(2D)$ the maps $W_{L,p_i}$ factor through
$\bigwedge^2 (H^0(C,L)/\lan 1\ran)$, since the section $1\in H^0(C,L)$ 
has a double zero at $p_i$. We have a decomposition
$$H^0(C,L)/\lan 1\ran=\bigoplus_{i=1}^g H^0(C,\OO(2p_i+D_i))/\lan 1\ran,$$
and an identification of each summand
$H^0(C,\OO(2p_i+D_i))/\lan 1\ran\simeq T_{p_i}^2$.
Now the restriction
$$T_{p_j}^2\ot T_{p_i}^2\hra \bigwedge^2 (H^0(C,L)/\lan 1\ran)\rTo{W_{L,p_i}} (L^2\ot\om_C)|_{p_i}
\simeq T_{p_i}^3$$
is given by $\a_{ji}$. Again, this determines $W_{L,p_i}$, since its restrictions
to $T_{p_j}^2\ot T_{p_k}^2$ are zero for  $j\neq i, k\neq i.$

\section{Reconstruction of the curve}\label{recon-sec}

In this section $(C,p_\bullet)$ corresponds to a generic point of $\MM_{g,g}$.
In particular, $h^0(D)=1$, where $D=p_1+\ldots+p_g$.

\subsection{Multiplication map}

Consider the line bundle 
$$L'=\OO_C(2D+p_1)=\OO_C(3p_1+2(p_2+\ldots+p_g))$$ 
of degree $2g+1$ on $C$. 

\begin{lem}\label{quadrics-lem} Let $g\ge 4$.
For generic $(C,p_\bullet)$ the curve $C$ is cut out by quadrics in the projective
embedding given by $|L'|$.
\end{lem}

\Pf . By \cite[Thm.\ 2]{GL}, this is true provided $C$ is not hyperelliptic and 
$L'\not\simeq \om_C(x+y+z)$ for any $x,y,z\in C$
(i.e., $C$ has no trisecants in the projective embedding given by $|L'|$) . Since $L'$ is determined by
$g\ge 4$ points on $C$, this holds generically.
\ed

Thus, for $g\ge 4$, 
generically we can recover the image of $C$ in $\P^{g+1}$ from the multiplication map
\begin{equation}\label{mult-map-eq}
H^0(C,L')\ot H^0(C,L')\to H^0(C,(L')^2).
\end{equation}

By the Riemann-Roch formula, we have $h^0(L')=g+2$, $h^0((L')^2)=3g+3$.

\begin{lem}\label{bases-lem}
For $i=1,\ldots,g$, let us pick a nonconstant rational function 
$f_i\in H^0(C,\OO(D+p_i))$ and a rational function
$h_i\in H^0(C,\OO(D+2p_i))\setminus H^0(C,\OO(D+p_i))$.
Then we have the following bases in $H^0(C,L')$ and $H^0(C,(L')^2)$:
\begin{align*}
& H^0(C,L'):\ 1, f_1, \ldots, f_g, h_1; \\
& H^0(C,(L')^2):\ 1, f_1, \ldots, f_g, h_1,\ldots, h_g, f_1^2,\ldots,f_g^2, f_1h_1, h_1^2.
\end{align*}
\end{lem}

\Pf . The exact sequences 
$$0\to H^0(C,\OO(nD))\to H^0(C,\OO((n+1)D))\to \bigoplus_{i=1}^g H^0(C,\OO((n+1)p_i)/\OO(np_i))
\to 0$$
for $n=1$, $2$ and $3$ give us the following bases:
\begin{align*}
& H^0(C,\OO(2D)):\ 1, f_1, \ldots, f_g; \\
& H^0(C,\OO(3D)):\ 1, f_1, \ldots, f_g, h_1,\ldots, h_g; \\
& H^0(C,\OO(4D)): 1, f_1, \ldots, f_g, h_1,\ldots, h_g, f_1^2,\ldots,f_g^2.
\end{align*}
Now the result follows from the exact sequences
$$0\to H^0(C,\OO(2D))\to H^0(C,L')\to H^0(C,\OO(3p_1)/\OO(2p_1))\to 0  \ \ \text{ and}$$
$$0\to H^0(C,\OO(4D))\to H^0(C,(L')^2)\to H^0(C,\OO(6p_1)/\OO(4p_1))\to 0.$$
\ed

We need convenient formal parameters at the marked points.

\begin{lem}\label{parameter-lem} 
Let $\cha(\k)=0$ (resp., $\cha(\k)>N$ for some $N$). 
Let $C$ be a smooth projective curve of genus $g$.
For any point $p$ and any divisor $E$ of degree $g-1$
such that $p\not\in\supp(E)$ and $h^0(p+E)=1$ there exists a formal parameter
$t_{p,E}$ (resp., formal parameter modulo $\mg^{N+1}$), 
unique up to rescaling by a constant, such that 
for every $n\ge 2$ (resp., for $2\le n\le N$),
there exists a global section of $\OO(np+E)$ with the polar part
$t_{p,E}^{-n}$ at $p$.
\end{lem}

\Pf . Pick a non-constant function $f(2)\in H^0(C,\OO(2p+E))$. 
Then for any local parameter $t$ at $p$ we can rescale $f(2)$ so that 
$$f(2)=\frac{1}{t^2}+ \frac{c}{t}+\ldots$$
at $p$, where $c$ depends only on $t\mod\mg^3$. Replacing $t$ by $t+at^2\mod\mg^3$
leads to the transformation $c\mapsto c-2a$. This implies the statement for $n=2$.
Then we proceed by induction: suppose we have a local parameter $t\mod\mg^n$ and functions
$f(m)\in H^0(C,\OO(mp+E))$ with polar parts $t^{-m}$ for $2\le m\le n-1$. 
Let us take $f(n)\in H^0(C,\OO(np+E))\setminus H^0(C,\OO((n-1)p+E))$.
Rescaling and subtracting an appropriate linear combination of $f(2),\ldots,f(n-1)$ we get
a unique such $f(n)$ with
$$f(n)=\frac{1}{t^n}+\frac{c}{t}+\ldots$$
at $p$, where we extend $t\mod\mg^n$ to $t\mod\mg^{n+1}$ in some way.
Changing $t$ by $t+at^n$ leads to the change of $c$ to $c-na$, so we find the unique
$t\mod\mg^{n+1}$ for which $c=0$.
\ed

Let $D_i=\sum_{j\neq i}p_i$.
Applying Lemma \ref{parameter-lem} (for $\cha(\k)\neq 2,3$) to the pairs $(p_i,D_i)$
we can choose formal parameters $(t_i=t_{p_i,D_i})$ at $p_i$, so that 
there are elements $f_i\in H^0(C,\OO(2p_i+D_i))$ and $h_i\in H^0(C,\OO(3p_i+D_i))$, such that
$$f_i\equiv \frac{1}{t_i^2}\mod\hat{\OO}_{C,p_i},$$
$$h_i\equiv\frac{1}{t_i^3}\mod\hat{\OO}_{C,p_i},$$
where $\hat{\OO}_{C,p_i}$ is the completion of the local ring $\OO_{C,p_i}$.
This fixes $f_i$ and $h_i$ up to adding a constant. 
Let $p_j$ be another marked point (so $j\neq i$). We have expansions
\begin{equation}
f_i\equiv\frac{\a_{ij}}{t_j}+\ga_{ij}+\de_{ij}t_j \mod t_j^2\hat{\OO}_{C,p_j},\nonumber
\end{equation}
\begin{equation}
h_i\equiv\frac{\b_{ij}}{t_j}+\varepsilon_{ij}+\vartheta_{ij}t_j \mod t_j^2\hat{\OO}_{C,p_j},\nonumber
\end{equation}
for some constants $(\a_{ij})$, $(\b_{ij})$, $(\ga_{ij})$, $(\de_{ij})$\footnote{We do not use 
the Kronecker delta in this paper}, $(\varepsilon_{ij})$, $(\vartheta_{ij})$ (defined for $i\neq j$).
Note that here $(\a_{ij})$ are the functions defined by \eqref{alpha-def-eq} (with some choices of trivializations
of the tangent spaces $T_{p_i}$).
Adding a constant to each $f_i$ (resp., $h_i$) we can
assume that
\begin{equation}
\ga_{i,i+1}=0, \ \ \varepsilon_{i,i+1}=0,
\end{equation}
for $i=1,\ldots,g$ (where we think of indices as elements of $\Z/g\Z$). This fixes the
choice of $f_i$ and $h_i$, for $i=1,\ldots,g$, uniquely.
Let us also define for each $i=1,\ldots,g$ a constant $\ga_{ii}$, so that at $p_i$ we have the expansion
\begin{equation}
f_i\equiv\frac{1}{t_i^2}+\ga_{ii}\mod t_i\hat{\OO}_{C,p_i}
\end{equation}
for some constants $(\ga_{ii})$.


To describe the multiplication map \eqref{mult-map-eq} in terms of the bases 
of Lemma \ref{bases-lem} we need to find
the decompositions of the products $f_if_j$ for $i\neq j$ and $f_ih_1$ for $i\neq 1$.

\begin{lem}\label{equations-lem} 
(i) 
For $i\neq j$ one has
\begin{equation}\label{ff-prod-eq}
f_if_j=\sum_{k\neq i,j}\a_{ik}\a_{jk}f_k+\a_{ji}h_i+\a_{ij}h_j+\ga_{ji}f_i+\ga_{ij}f_j+a_{ij},
\end{equation}
for some constants $a_{ij}=a_{ji}$. Furthermore, one has the following relations: 
\begin{equation}\label{bc-rel}
\de_{ij}+\a_{ij}\ga_{jj}=\sum_{k\neq i,j}\a_{ik}\a_{jk}\a_{kj}+\a_{ji}\b_{ij}+\ga_{ji}\a_{ij}, 
\end{equation}
\begin{equation}\label{ct-rel}
\a_{ik}(\ga_{jk}-\ga_{ji})+\a_{jk}(\ga_{ik}-\ga_{ij})=\sum_{l\neq i,j,k}\a_{il}\a_{jl}\a_{lk}+\a_{ji}\b_{ik}+\a_{ij}\b_{jk},
\end{equation}
\begin{equation}\label{abcd-rel}
\a_{ik}\de_{jk}+\a_{jk}\de_{ik}+\ga_{ik}\ga_{jk}=
\sum_{l\neq i,j}\a_{il}\a_{jl}\ga_{lk}+\a_{ji}\varepsilon_{ik}+\a_{ij}\varepsilon_{jk}+\ga_{ji}\ga_{ik}+\ga_{ij}\ga_{jk}+a_{ij},
\end{equation}
where $i,j,k$ are distinct (but $l=k$ is allowed in the last equation).

\noindent
(ii) For $i\neq j$ one has
\begin{equation}\label{fg-prod-eq}
f_ih_j=\sum_{k\neq i,j} \a_{ik}\b_{jk}f_k+\b_{ji}h_i+\ga_{ij}h_j+\varepsilon_{ji}f_i+\de_{ij}f_j+
\a_{ij}[f_j^2-2\ga_{jj}f_j-\sum_{k\neq j}\a_{jk}^2f_k]+b_{ij}
\end{equation}
for some constants $(b_{ij})$. Furthermore, one has the following relations
\begin{equation}\label{pr-rel}
\vartheta_{ji}+\b_{ji}\ga_{ii}=\sum_{k\neq i,j}\a_{ik}\b_{jk}\a_{ki}+\ga_{ij}\b_{ji}+\de_{ij}\a_{ji}
+\a_{ij}[2\a_{ji}\ga_{ji}-2\ga_{jj}\a_{ji}-\sum_{k\neq i,j}\a_{jk}^2\a_{ki}],
\end{equation}
\begin{equation}\label{r-rel}
\begin{array}{l}
\a_{ik}\varepsilon_{jk}+\ga_{ik}\b_{jk}=\sum_{l\neq i,j,k}\a_{il}\b_{jl}\a_{lk}+\b_{ji}\b_{ik}+
\ga_{ij}\b_{jk}+\varepsilon_{ji}\a_{ik}+\de_{ij}\a_{jk}+\\
\a_{ij}[2\a_{jk}\ga_{jk}-2\ga_{jj}\a_{jk}-\sum_{l\neq j,k}\a_{jl}^2\a_{lk}],
\end{array}
\end{equation}
\begin{equation}\label{pqe-rel}
\begin{array}{l}
\a_{ik}\vartheta_{jk}+\ga_{ik}\varepsilon_{jk}+\de_{ik}\b_{jk}=
\sum_{l\neq i,j}\a_{il}\b_{jl}\ga_{lk}+\b_{ji}\varepsilon_{ik}+\ga_{ij}\varepsilon_{jk}+\varepsilon_{ji}\ga_{ik}+\de_{ij}\ga_{jk}+\\ 
\a_{ij}[\ga_{jk}^2+2\a_{jk}\de_{jk}-2\ga_{jj}\ga_{jk}-\sum_{l\neq j}\a_{jl}^2\ga_{lk}]
+b_{ij},
\end{array}
\end{equation}
where $i,j,k$ are distinct.
\end{lem}

\Pf . (i) We have $f_if_j\in H^0(C,\OO(3p_i+3p_j+2D_{ij}))$. Expanding in the formal parameter at
$p_i$ we obtain
\begin{equation}
f_if_j= (\frac{1}{t_i^2}+\ga_{ii}+\ldots)(\frac{\a_{ji}}{t_i}+\varepsilon_{ji}+\de_{ji}t+\ldots)=
\frac{\a_{ji}}{t_i^3}+\frac{\ga_{ji}}{t_i^2}+\frac{\de_{ji}+\a_{ji}\ga_{ii}}{t_i}+\ldots.
\end{equation}
Hence, the difference
$$f_if_j-\a_{ji}h_i-\a_{ij}h_j-\ga_{ji}f_i-\ga_{ij}f_j$$
has poles of order at most $1$ at $p_i$ and $p_j$. On the other hand, expanding at $p_k$, where 
$k\neq i,j$ we obtain
\begin{equation}
f_if_j=\frac{\a_{ik}\a_{jk}}{t_k^2}+\frac{\a_{ik}\ga_{jk}+\a_{jk}\ga_{ik}}{t_k}+
\a_{ik}\de_{jk}+\a_{jk}\de_{ik}+\ga_{ik}\ga_{kj} \mod t_k\hat{\OO}_{C,p_k}.
\end{equation}
Hence,
$$f_if_j-\sum_{k\neq i,j}\a_{ik}\a_{jk}f_k-\a_{ji}h_i-\a_{ij}h_j-\ga_{ji}f_i-\ga_{ij}f_j$$
has poles of order at most $1$ at all marked points. Since such a function has to be
constant, this implies \eqref{ff-prod-eq}. Now the relation \eqref{bc-rel} is obtained
by equating polar parts of both sides of \eqref{ff-prod-eq} at $p_j$. Similarly,
\eqref{ct-rel} and \eqref{abcd-rel} are obtained by considering expansions of both sides
of \eqref{ff-prod-eq} at $p_k$, where $k\neq i,j$.

\noindent
(ii) To prove \eqref{fg-prod-eq} we start by observing
that 
$$f_j^2-2\ga_{jj}f_j-\sum_{k\neq j}\a_{jk}^2f_k\equiv \frac{1}{t_j^4}\mod t_j^{-1}\hat{\OO}_{C,p_j}$$
and that this function has poles of order at most $1$ at all $p_k$ for $k\neq j$. Then we argue as in the proof
of \eqref{fg-prod-eq}.
Comparing the polar parts of both sides of \eqref{fg-prod-eq} at $p_i$ we obtain
\eqref{pr-rel}. The relations \eqref{r-rel} and \eqref{pqe-rel} are obtained by considering
expansions of both sides of \eqref{fg-prod-eq} at $p_k$, where $k\neq i,j$.
\ed

\noindent 
{\it Proof of Theorem \ref{birational-thm}.} 
We would like to prove that the map $\a:\MM^{(1)}_{g,g}\rTo{(\a_{ij})}\G_m^{g^2-g}$ is generically one-to-one
on its image for $g\ge 6$. Since the restriction of $\a$ to
fibers of the projection $\MM^{(1)}_{g,g}\to\MM_{g,g}$ is injective, it is enough to show how to recover 
generic $(C,p_\bullet)$ from the constants $(\a_{ij})$, defined using some trivializations of the tangent spaces
$T_{p_i}$.
By Lemma \ref{quadrics-lem}, for generic $(C,p_\bullet)$ the kernel of the multiplication map
\eqref{mult-map-eq} gives quadratic equations which cut out $C$ in the projective
embedding given by $|2D+p_1|$, where $D=p_1+\ldots+p_g$. 
If in addition we know the line spanned by the section $1\in H^0(C,\OO(2D+p_1))$
then we can recover $p_1$ as a triple zero of this section and the unordered collection of 
points $p_2,\ldots,p_g$ as double zeros of this section. Furthermore,
we can recover each $p_i$ for $i\ge 2$ if we know the line spanned by the section
$f_i\in H^0(C,\OO(D+p_i))\sub H^0(C,\OO(2D+p_1))$ used in Lemma \ref{bases-lem}.
Indeed,  generically $f_i$, viewed as a section of $L'$, is nonzero near $p_i$ and has simple zeros
at $p_j$ for $j\neq i$, $j\ge 2$.
 
By \eqref{ff-prod-eq} and \eqref{fg-prod-eq},
the constants $(\a_{ij})$, $(\b_{ij})$, $(\ga_{ij})$, $(\de_{ij})$, 
$(\varepsilon_{ij})$, $(a_{ij})$ and $(b_{ij})$ (where $i\neq j$)
determine the multiplication map \eqref{mult-map-eq} with respect to the bases of Lemma
\ref{bases-lem}. 
Thus, it is enough to show that for generic 
$(C,p_\bullet)\in\ov{\MM}_{g,g}$ 
these constants are uniquely determined by
$(\a_{ij})$. We do this by solving the equations obtained in Lemma \ref{equations-lem}.

\noindent
{\bf Step 1}. We would like to solve the equations \eqref{ct-rel} (together with
the condition $\ga_{i,i+1}=0$) for $(\b_{ij})$, $(\ga_{ij})$. The fact that for generic
$(C,p_\bullet)$ these equations
determine $(\b_{ij})$ and $(\ga_{ij})$ follows from Proposition \ref{nondeg-prop}(i) below.

\noindent
{\bf Step 2}. Note that we can express $\de_{ij}$ in terms of $\ga_{jj}$ (and known quantities) using \eqref{bc-rel}.
Now substituting in \eqref{abcd-rel} 
we get a linear system for $(\varepsilon_{ij})$, $(a_{ij})$ and $(\ga_{ii})$ of
the form
\begin{equation}\label{acd-rel}
\a_{ji}\varepsilon_{ik}+\a_{ij}\varepsilon_{jk}+3\a_{ik}\a_{jk}\ga_{kk}+a_{ij}=A_{ijk}
\end{equation}
By Proposition \ref{nondeg-prop}(ii) below, for generic $(C,p_\bullet)$
these equations (together with the condition $\varepsilon_{i,i+1}=0$) determine
$(\varepsilon_{ij})$, $(a_{ij})$ and $(\ga_{ii})$ uniquely.

\noindent
{\bf Step 3}. Using \eqref{pr-rel} we can express $\vartheta_{ij}$ in terms of the constants determined
by $\a_{ij}$ (generically). Finally, we use \eqref{pqe-rel} to express $(b_{ij})$ in terms of the known
constants.
\ed

\begin{rem} The above reconstruction procedure gives also a way to
produce some polynomial equations for $(\a_{ij})$ for $g\ge 6$. For example, 
for $g=6$ the system \eqref{ct-rel} of $60$ equations has $54$ variables
(since we set $\ga_{i,i+1}=0$). Taking any $55$ equations we get the vanishing
of a $55\times 55$ determinant, with one column of homogeneous cubic polynomials
in $(\a_{ij})$ and all other entries linear in $(\a_{ij})$, which gives a degree-$57$ equation
on $(\a_{ij})$. One can check that this indeed leads to nonzero equations. 
\end{rem}

\subsection{Degeneration argument}\label{degen-sec}

So far, we have reduced our reconstruction problem to proving that certain linear systems
have maximal possible rank generically on $\ov{\MM}_{g,g}$. Namely, consider
the homogeneous linear system on $(\b_{ij}, \ga_{ij})$, associated with \eqref{ct-rel},
\begin{equation}\label{ct-hom-rel}
\a_{ik}(\ga_{jk}-\ga_{ji})+\a_{jk}(\ga_{ik}-\ga_{ij})=\a_{ji}\b_{ik}+\a_{ij}\b_{jk},
\end{equation} 
and the homogeneous linear system on $(\varepsilon_{ij}, a_{ij}, \ga_{ii})$, associated with \eqref{acd-rel},
\begin{equation}\label{acd-hom-rel}
\a_{ji}\varepsilon_{ik}+\a_{ij}\varepsilon_{jk}+3\a_{ik}\a_{jk}\ga_{kk}+a_{ij}=0,
\end{equation}
(in both systems $i,j,k$ are distinct).
We have to check that generically they have only the obvious solutions
\begin{equation}\label{ct-triv-sol}
\ga_{ij}=\la_i, \ \ \b_{ij}=0,
\end{equation}
\begin{equation}\label{acd-triv-sol}
\varepsilon_{ij}=-\mu_i, \ \ a_{ij}=\a_{ji}\mu_i+\a_{ij}\mu_j, \ \ \ga_{ii}=0,
\end{equation}
for some $(\la_i)$ and $(\mu_i)$. Our strategy is to reduce this to the case $g=6$ and to
study the above systems for irreducible rational nodal curves, for which $\a_{ij}$ can be determined
explicitly.

Namely, consider the curve $C$ obtained from $\P^1$ by gluing $g$ pairs of distinct points 
$(a_1,b_1),\ldots,(a_g,b_g)$, where $a_i,b_i\in\A^1$, together with the marked points 
$p_1,\ldots,p_g\in C$ that are images of the points $c_1,\ldots,c_g\in\A^1\sub\P^1$. 
Note that the coordinate on $\A^1$ gives rise
to a trivialization of the tangent line to $C$ at each $p_i$.
We look for the rational functions $f_i\in H^0(C,\OO(2p_i+D_i))$ in the form
$$f_i(t)=\frac{1}{(t-c_i)^2}+\sum_{j=1}^g \frac{\a_{ij}}{t-c_j},$$
where for $i\neq j$ the constants $\a_{ij}$ are the functions we are interested in 
(while $\a_{ii}$ do not have an invariant meaning). Note that to compute $\a_{ij}$
we do not need the special choice of parameters at $p_i$ that we used for Lemma \ref{equations-lem}.
The conditions $f_i(a_k)=f_i(b_k)$ give the following system of linear equations
on $(\a_{ij})$:
$$\sum_{j=1}^g\left(\frac{1}{b_k-c_j}-\frac{1}{a_k-c_j}\right)\a_{ij}=
\frac{1}{(a_k-c_i)^2}-\frac{1}{(b_k-c_i)^2}, \ 1\le i,k\le g.$$
Dividing by $a_k-b_k$ we can rewrite this as
\begin{equation}\label{alpha-rel}
\sum_{j=1}^g\frac{1}{(b_k-c_j)(a_k-c_j)}\a_{ij}=
\frac{2c_i-a_k-b_k}{(a_k-c_i)^2(b_k-c_i)^2}, \ 1\le i,k\le g.
\end{equation}
Let us consider the $g\times g$-matrices 
$A=(\a_{ij})$, $M=(m_{ij})$ and $N=(n_{ij})$, where
\begin{equation}\label{m-n-eq}
m_{ij}=\frac{1}{(b_j-c_i)(a_j-c_i)},\ \ n_{ij}=
\frac{2c_i-a_j-b_j}{(a_j-c_i)^2(b_j-c_i)^2}.
\end{equation}
Then \eqref{alpha-rel} is simply the matrix relation
$$AM=N.$$
Note that the matrices $M$ and $N$ are also defined when $a_i=b_i$.

\begin{lem}\label{comp-lem} Let $\cha(\k)=0$.
For $g=6$ and
$$a_i=b_i=-c_i=i \text{ for } i=1,\ldots,6,$$
the matrix $M$ is invertible. Furthermore, for the corresponding
entries $\a_{ij}$ of the matrix
$A=NM^{-1}$ each of the systems \eqref{ct-hom-rel} and \eqref{acd-hom-rel} has
$6$ free variables. Hence, the same assertion is true for generic $a_i,b_i,c_i$.
\end{lem}

\Pf . We checked this with the help of a computer (see Appendix).
\ed

\begin{prop}\label{nondeg-prop} Let $\cha(\k)=0$ and $g\ge 6$.

\noindent
(i) At generic point of $\ov{\MM}^{(1)}_{g,g}$ 
the system \eqref{ct-hom-rel} has only trivial solutions \eqref{ct-triv-sol}.

\noindent
(ii) At generic point of $\ov{\MM}^{(1)}_{g,g}$ 
the system \eqref{acd-hom-rel} has only trivial solutions \eqref{acd-triv-sol}.
\end{prop}

\Pf . (i) Lemma \ref{comp-lem} implies that the assertion is true
for generic $(C,p_1,\ldots,p_6)\in\MM_{6,6}^{(1)}$.

For $g>6$ let us fix a subset $I\sub\{1,\ldots,g\}$ consisting of $6$ elements. 
We claim that generically the only solution of
the equations \eqref{ct-hom-rel} with $i,j,k\in I$ for variables $(\b_{ij},\ga_{ij}\ | i,j\in I)$ is 
\begin{equation}\label{tc-sol}
\b_{ij}=0, \ \ \ga_{ij}=\la_{i,I},
\end{equation}
for some constants $\la_{i,I}$.
Indeed, without loss of generality we can assume that $I=\{1,\ldots,6\}$.
Let us take generic curves $(C_1,p_1,\ldots,p_6)\in\MM_{6,6}$ and 
$(C_2,p_7,\ldots,p_g)\in\MM_{g-6,g-6}$ and 
consider the nodal curve $(C,p_1,\ldots,p_g)$ obtained from $C_1\sqcup C_2$
by identifying points $p\in C_1$ and $q\in C_2$
(where $p$ and $q$ are different from all the markings). We also assume that nonzero
tangent vector fields are chosen at all points, so $\a_{ij}$ are defined for $i\neq j$.
Now for $i\in\{1,\ldots,6\}$ a nonconstant section $f_i\in H^0(C,\OO(2p_i+D_i))$ will
restrict to a similar section on $C_1$ (and will have a constant restriction to $C_2$).
Hence for $i,j\in \{1,\ldots,6\}$ the constants $\a_{ij}$ calculated for $(C,p_\bullet)$
are equal to those for $(C_1,p_1,\ldots,p_6)$. This implies our claim.

Thus, for generic $(C,p_\bullet)$ a solution of \eqref{ct-rel} satisfies \eqref{tc-sol} for
each $I$ as above.
Note that for $I$ and $I'$ such that $i,j\in I\cap I'$ we have $\la_{i,I}=\la_{i,I'}=\ga_{ij}$.
Since any two subsets $I$ containing $i$ can be connected by a chain of subsets containing $i$,
in which every two consecutive terms have at least two elements in common, this implies
that $\la_{i,I}$ depends only on $i$.

\noindent
(ii) For $g=6$ this follows from Lemma \ref{comp-lem}. Then we proceed as in part (i).
\ed

\section{The tangent map}\label{tangent-map-sec}

\subsection{General formula}

The tangent space to the moduli space $\ov{\MM}^{(1)}_{g,g}$ 
at a stable curve $(C,p_1,\ldots,p_g, v_1,\ldots,v_g)$ (where $v_i\in T_{p_i}\setminus 0$)
is canonically identified with
$\Ext^1(\Om_C,\OO(-2D))$, where $\Om_C$ denotes the sheaf of K\"ahler differentials and $D=\sum_{i=1}^g p_i$.
On the other hand, if $h^0(\OO(D))=1$ then we can
use the boundary homomorphism
$$\bigoplus_{j\neq i} T_{p_j}\simeq H^0(C,\OO(D)/\OO(p_i))\rTo{\sim} H^1(C,\OO(p_i))$$
to get natural bases in each space $H^1(C,\OO(p_i))$ numbered by $e_{ij}$, $i\neq j$.

Let us consider the regular map
$$\a^{reg}:\UU^{(1)}\rTo{(\a_{ij})}\A^{g^2-g},$$
where $\UU^{(1)}\sub\ov{\MM}^{(1)}_{g,g}$ is the preimage of the open substack
$\UU\sub\ov{\MM}_{g,g}$ defined by $h^0(D)=1$.

\begin{prop}\label{tangent-prop}
Under the above identifications the tangent map to $\a^{reg}$
is the map
\begin{equation}\label{tangent-map-eq}
\Ext^1(\Om_C,\OO(-2D))\rTo{(df_i)} \bigoplus_{i=1}^g\Ext^1(\OO(-2D-p_i),\OO(-2D))\simeq
\bigoplus_{i=1}^g H^1(C,\OO(p_i)),
\end{equation}
where for each $i$, $df_i\in\Om_C(2D+p_i)$ is the differential of the rational function 
$f_i\in H^0(\OO(D+p_i))$, such that $f_i\equiv v_i^2\mod\OO(D)$.
\end{prop}

\Pf . By irreducibility of the moduli space it is enough to consider the case when
$C$ is smooth. Then \eqref{tangent-map-eq} is the map induced on $H^1$
by the morphism of coherent sheaves
$$\TT(-2D)\rTo{(df_i)}\bigoplus_i \OO(p_i).$$
Note also that for each $i$ the natural map
$$\OO(p_i)\to\bigoplus_{j\neq i}\OO(D_j)$$
induces an isomorphism on $H^1$ (recall that $D_j=D-p_j$).
Hence, our assertion reduces to checking that for each $i\neq j$,
the differential $d\a_{ij}$ of the function $\a_{ij}$ at $(C,p_1,\ldots,p_g,v_1,\ldots,v_g)$
is equal to the map induced on $H^1$ by the morphism
\begin{equation}\label{T-2D-Dj-eq}
\TT(-2D)\rTo{df_i}\OO(p_i)\to\OO(D_j).
\end{equation}
To this end let us fix an affine covering $(U_a)$ of $C$ and a Cech $1$-cocycle
$v_{ab}$ with values in $\TT(-2D)$ (we 
assume that each marked point is contained in only one $U_a$). 
This gives a first-order deformation of $(C,p_1,\ldots,p_g)$
over $\k[\eps]/(\eps^2)$, glued from the trivial deformations 
$U_a[\eps]:=U_a\times\Spec(\k[\eps]/(\eps^2))$ of $U_a$ with the
transitions on $U_{ab}[\eps]$ given by the automorphisms $\id+\eps v_{ab}$.
Let $f_i\in H^0(C,\OO(D+p_i))$ be such that $f_i\equiv v_i^2\mod\OO(D)$, so that
by definition
$$f_i\equiv \a_{ij}\cdot v_j\mod \OO(D_j+p_i).$$
We want to deform this function over $\k[\eps]/(\eps^2)$ preserving the condition
$f_i\equiv v_i^2\mod\OO(D)$. Thus, the deformed function should have form
$f_i+\eps g_a$ on each $U_a$, where $g_a\in H^0(U_a,\OO(D))$.
The gluing condition gives
$$f_i+\eps g_b=(\id+\eps v_{ab})(f_i+\eps g_a)$$
on $U_{ab}[\eps]$, which is equivalent to
\begin{equation}\label{g-ab-eq}
g_b=g_a+v_{ab}(f_i).
\end{equation}
Then $d\a_{ij}$ is the image of $g_{a(j)}$ in $\OO(p_j)/\OO$, where $p_j\in U_{a(j)}$. 
From the exact sequence 
$$0\to \OO(D_j)\to \OO(D)\to \OO(p_j)/\OO\to 0$$
we get the following exact sequence of Cech complexes:
\begin{diagram}
0&\rTo{} &C^0(\OO(D_j))&\rTo{}& C^0(\OO(D))&\rTo{r_j}&C^0(\OO(p_j)/\OO)&\rTo{} 0\\
&&\dTo{}&&\dTo{\de}&&\dTo{}\\
0&\rTo{} &C^1(\OO(D_j))&\rTo{\iota}& C^1(\OO(D))&\rTo{}& 0 &\rTo{} 0
\end{diagram}
Note that 
$$v_{ab}(f_i)=\lan v_{ab}, df_i\ran\in \OO(p_i)\sub\OO(D_j),$$
so we can view $(v_{ab}(f_i))$ as a $1$-cocycle in $C^1(\OO(D_j))$.
By \eqref{g-ab-eq}, the cochain $(g_a)\in C^0(\OO(D))$ satisfies
$$\de((g_a))=\iota(v_{ab}(f_i)).$$
Since $r_j((g_a))$ is exactly $d\a_{ij}$, we obtain that the class
$[v_{ab}(f_i)]\in H^1(C,\OO(D_j))$ is the image of $d\a_{ij}$ under the connecting
homomorphism $H^0(\OO(p_j)/\OO)\to H^1(\OO(D_j))$.
Since the map \eqref{T-2D-Dj-eq} is given by $v\mapsto v(f_i)$, this implies our claim.
\ed

\subsection{Tangent map at a rational irreducible nodal curve}

Let $(C,p_1,\ldots,p_g)$ be a stable curve.
Using Serre duality we can identify the dual to the tangent map \eqref{tangent-map-eq}
with the map
\begin{equation}\label{dual-tangent-map-eq}
\bigoplus_{i=1}^g H^0(C,\om_C(-p_i))\rTo{(df_i)} H^0(C,\Om_C\ot\om_C(2D)),
\end{equation}
where $\om_C$ is the dualizing sheaf on $C$.

We are going to describe the map \eqref{dual-tangent-map-eq} explicitly in the
case of the curve $C$ obtained from $\P^1$ by gluing $g$ pairs of distinct points
$(a_1,b_1),\ldots,(a_g,b_g)$, with the marked points $c_1,\ldots,c_g$.
Let us  denote by $q_i\in C$ the node corresponding to a pair $(a_i,b_i)$.
Recall that in this case the functions $f_i\in\OO(2p_i+D_i)$ correspond to the functions on
$\P^1$
$$f_i(t)=\frac{1}{(t-c_i)^2}+\sum_{j=1}^g \frac{\a_{ij}}{t-c_j},$$
where the matrix $(\a_{ij})$ is determined by the conditions $f_i(a_k)=f_i(b_k)$ 
(see Section \ref{degen-sec}).

The main problem is to understand the space $H^0(C,\Om_C\ot\om_C(2D))$.
Recall (see \cite[Exer.\ 5.9]{Ha-DT}) that $\Om_C$ fits into the exact sequence 
$$0\to\OO_Z\to\Om_C\rTo{\nu}\om_C\to\OO_Z\to 0,$$
where $Z$ is the union of all nodes. Formally at a node $q\in C$ the curve looks as
$\Spec(\k[[x,y]]/(xy))$ and the completion of $\Om_C$ is generated by $dx$ and $dy$ with the 
relation
$xdy=-ydx$, so that the embedding $\OO_Z\to\Om_C$ is given by $1\mapsto ydx$.
The dualizing sheaf $\om_C$ is locally free and is generated near the node by $dx/x=-dy/y$. 
Thus, the space $(\Om_C\ot\om_C)\ot\OO_{C,q}/\mg_q^2$ 
has the basis $dx\ot (dx/x), dy\ot (dx/x), ydx\ot (dx/x), xdx\ot (dx/x), ydy\ot (dx/x)$.
Let
$$\tau_{q}:(\Om_C\ot\om_C)\ot\OO_{C,q}/\mg_q^2\to \k\cdot (xdy\ot\frac{dx}{x})$$
denote the projection to $(xdy \ot (dx/x))$ with respect to this basis.
Then we have an embedding
\begin{equation}\label{nodal-emb-1-eq}
H^0(C,\Om_C\ot\om_C(2D))\rTo{\nu, (\tau_q)} H^0(C,\om_C^{\ot 2}(2D))\oplus \k^{Z}.
\end{equation}
Thus, to describe \eqref{dual-tangent-map-eq} it is enough to consider its compositions
with $\nu$ and with $\tau_q$ for all $q\in Z$.

Let $\pi:\P^1\to C$ be the normalization map. Then we have an isomorphism
$$\pi^*\om_C\simeq\om_{\P^1}(\sum_i (a_i+b_i)),$$
so that near the $i$th node $q_i$ the sections of $\om_C$ are distinguished by the
condition $\Res_{a_i}+\Res_{b_i}=0$.
In particular, we have an embedding
$$H^0(C,\om_C^{\ot 2}(2D))\hra H^0(\P^1,\om_{\P^1}^{\ot 2}(2D+2\sum_i (a_i+b_i))).$$
Since the pull-backs $\pi^*df_i\in\om_{\P^1}(p_i+2D)$ are regular at all $a_i$'s and $b_i$'s,
we have a commutative diagram
\begin{diagram}
\bigoplus_{i=1}^g H^0(C,\om_C(-p_i))&\rTo{\nu\circ(df_i)}&
H^0(C,\om_C^{\ot 2}(2D))\\
\dTo{\varphi}&&\dTo{}\\
H^0(\P^1,\om_{\P^1}^{\ot 2}(2D+\sum_j (a_j+b_j)))&\rTo{}
&H^0(\P^1,\om_{\P^1}^{\ot 2}(2D+2\sum_j (a_j+b_j)))
\end{diagram}
where $\varphi$ is induced by the embeddings
$$H^0(C,\om_C(-p_i))\hra H^0(\P^1,\om_\P^1(-c_i+\sum_j (a_j+b_j)))$$
followed by the product with $\pi^*df_i$.
Finally, since $\om_{\P^1}^{\ot 2}$ has no global sections, we have an embedding 
$$\iota:H^0(\P^1,\om_{\P^1}^{\ot 2}(2D+\sum_j (a_j+b_j)))\to\bigoplus_{i=1}^g\om_{\P^1}^{\ot 2}(2c_i)/\om^{\ot 2}_{\P^1}
\oplus\bigoplus_{j=1}^g\om^{\ot 2}_{\P^1}(a_j)/\om^{\ot 2}_{\P^1}\oplus\om^{\ot 2}_{\P^1}(b_j)/\om_{\P^1},
$$
given by the polar parts at all points $a_i$, $b_i$ and $c_i$.
Thus, the map \eqref{dual-tangent-map-eq} is essentially determined by the map
\begin{equation}\label{nodal-dual-tangent-eq}
\bigoplus_{i=1}^g H^0(C,\om_C(-p_i))\rTo{\iota\circ\varphi,(\tau_q)} 
\bigoplus_{i=1}^g\om_{\P^1}^{\ot 2}(2c_i)/\om^{\ot 2}_{\P^1}
\oplus\bigoplus_{j=1}^g\om^{\ot 2}_{\P^1}(a_j)/\om^{\ot 2}_{\P^1}\oplus\om^{\ot 2}_{\P^1}(b_j)/\om_{\P^1}\oplus \k^Z.
\end{equation}
In particular, \eqref{nodal-dual-tangent-eq} has the same rank as
\eqref{dual-tangent-map-eq}, which is the same as the rank of the tangent map \eqref{tangent-map-eq}.

Now we proceed to calculating the components of \eqref{nodal-dual-tangent-eq} explicitly.
First, let us describe a basis in the source space. We realize global sections of $\om_C$
as rational $1$-forms on $\P^1$ having $1$st order poles at $a_i$ and $b_i$ with opposite residues
at $a_i$ and $b_i$ for each $i$ (and regular elsewhere). All such $1$-forms on $\P^1$ can be 
written as
$$\sum_{k=1}^g \frac{x_k}{(t-a_k)(t-b_k)}dt,$$
for some constants $x_1,\ldots,x_g$. Thus if $(x_{ij})=M^{-1}$, where $M=(m_{ij})$ 
(see \eqref{m-n-eq}),
then for each $j$ the form
\begin{equation}\label{eta-e-eq}
\eta_j=e_j(t)dt, \text{ where } e_j(t)=\sum_{k=1}^g \frac{x_{kj}}{(t-a_k)(t-b_k)},
\end{equation}
is a generator of the $1$-dimensional subspace $H^0(C,\om_C(-D_j))\sub H^0(C,\om_C)$.
Thus, we can take $(\eta_j)_{j\neq i}$ as a basis of $H^0(C,\om_C(-p_i))$.
It remains to calculate the polar parts of $df_i\ot\eta_j$, where $i\neq j$, at all the points
$a_k$, $b_k$ and $c_k$, as well as the constants
$\tau_{q_k}(df_i\ot\eta_j)$. The latter constants can computed as follows.

\begin{lem} Let $U$ be a neighborhood of a node $q\in C$, and let
$x$ and $y$ be formal parameters at the two points $a$ and $b$ over the node on the normalization. 
For $\eta\in H^0(U,\om_U)$ consider the expansions near $a$ and $b$ of
its pull-back to the normalization,
$$\wt{\eta}=(d_{-1}+d_0x+\ldots)\frac{dx}{x}, \ \ \wt{\eta}=(e_{-1}+e_0y+\ldots)\frac{dy}{y},$$
where $d_{-1}+e_{-1}=0$.
Then for $f\in\OO(U)$ we have
$$\tau_q(df\ot\eta)=\left(e_0\frac{d\wt{f}}{dx}(a)+d_0\frac{d\wt{f}}{dy}(b)\right)\cdot xdy\ot\frac{dx}{x},$$
where $\wt{f}$ is the pull-back of $f$ to the normalization.
\end{lem}

\Pf . Let $\wt{f}=P(x)$ at $a$ and $\wt{f}=Q(y)$ at $b$, where $P\in \k[[x]]$, $Q\in \k[y]]$,
$P(0)=Q(0)$. Then 
$$df=P'(x)dx+Q'(y)dy\in\Om_{C,q}\ot\hat{\OO}_{C,q}.$$
Under the trivialization of $\om_C$ in the formal neighborhood $q$ given by $dx/x$,
$\eta$ corresponds to 
$$(d_{-1}+d_0x+d_1x+\ldots-e_0y-e_1y-\ldots)\frac{dx}{x}.$$
Hence
$$df\ot\eta=(d_{-1}+d_0x-e_0y+\ldots)(P'(x)dx+Q'(y)dy)\ot\frac{dx}{x}.$$
The terms contributing to $\tau_q$ are 
$$(d_0 Q'(0)xdy-e_0P'(0)ydx)\ot\frac{dx}{x}=(d_0Q'(0)+e_0P'(0))ydx\ot\frac{dx}{x}$$
which gives the result.
\ed

To apply this lemma in our case we use
expansions of $\eta_j$ near $a_k$ and $b_k$:
$$\eta_j=\left(\frac{x_{kj}}{(a_k-b_k)(t-a_k)}+[\sum_{l\neq k}\frac{x_{lj}}{(a_k-a_l)(a_k-b_l)}-
\frac{x_{kj}}{(a_k-b_k)^2}]+\ldots\right)dt,$$
$$\eta_j=\left(\frac{x_{kj}}{(b_k-a_k)(t-b_k)}+[\sum_{l\neq k}\frac{x_{lj}}{(b_k-a_l)(b_k-b_l)}-
\frac{x_{kj}}{(a_k-b_k)^2}]+\ldots\right)dt.$$
Hence,
$$\tau_{q_k}(df_i\ot\eta_j)=\wt{e}_{jk}(a_k)\cdot f'_i(b_k)+\wt{e}_{jk}(b_k)\cdot f'_i(a_k),$$
where
$$\wt{e}_{jk}(t)=\sum_{l\neq k}\frac{x_{lj}}{(t-a_l)(t-b_l)}-\frac{x_{kj}}{(a_k-b_k)^2}.$$
Calculation of the polar parts is straightforward. The polar part of $df_i\ot\eta_j$ at
$a_k$ (resp., $b_k$) is
$$\frac{x_{kj}f'_i(a_k)}{a_k-b_k}\cdot \frac{dt^{\ot 2}}{t-a_k}
\ \ \left(\text{resp.,}\  \frac{x_{kj}f'_i(a_k)}{a_k-b_k}\cdot \frac{dt^{\ot 2}}{t-a_k}\right).$$
To calculate polar parts at $c_k$ we need expansions of $f_i$ and $\eta_j$ in $t-c_k$,
so these will be expressed in terms of $\a_{ik}$ and of first two
derivatives of $e_j(t)$ at $c_k$ (see \eqref{eta-e-eq}). Namely for $k\neq i,j$ the polar part
of $df_i\ot\eta_j$ at $c_k$ is
$$\frac{e'_j(c_k)\a_{ik}}{t-c_k} (dt)^{\ot 2}$$
The polar part of $df_i\ot\eta_j$ at $c_i$ is
$$\left(\frac{2e'_j(c_i)}{(t-c_i)^2}+\frac{e''_j(c_i)+e'_j(c_i)\a_{ii}}{t-c_i}\right) dt^{\ot 2}.$$
Finally, the polar part of $df_i\ot\eta_j$ at $c_j$ is given by
$$\left(\frac{\a_{ij}}{(t-c_j)^2}+\frac{e'_j(c_j)\a_{ij}}{t-c_j}\right) dt^{\ot 2}.$$

\begin{thm}\label{dominant-thm} Assume that $\cha(\k)=0$. 
For $g\le 5$ the rational map 
$\ov{\a}:\MM_{g,g}\to\G_m^{g^2-2g}$
is dominant.
\end{thm}

\Pf . For $g=3$ this follows from Proposition \ref{image-prop}, which in this case states
that the restriction of $\ov{\a}$ to the generic fiber of the projection $\MM_{g,g}\to\MM_g$ is
generically \'etale. For $g=4$ and $g=5$ we use the above calculation to
construct a rational irreducible nodal curve with $g$ points for which the tangent map
to $\a$ (and hence to $\ov{\a}$) is surjective. Namely, we check using the computer that for
$a_i=-c_i=i$, $b_i=g+i$, where $g=4$ or $5$, the rank of the map
\eqref{nodal-dual-tangent-eq} is $g^2-g$ (see the GAP codes in the Appendix),
hence the tangent map \eqref{tangent-map-eq} has the same rank.
\ed 

\begin{rem} By Proposition \ref{image-prop}, in the case $g=3$ the map $\ov{\a}$ is still
dominant and generically smooth for $\cha(\k)>0$. In the cases $g=4$ and $g=5$ the
same is true if $\cha(\k)$ is sufficiently large.
\end{rem} 

\bigskip

\centerline{APPENDIX. GAP codes}

\bigskip

\noindent

{\bf 1. GAP codes for Lemma \ref{comp-lem}}.

\medskip

Setting up vectors $a=(a_i)$, $b=(b_i)$, $c=(c_i)$
and calculating the matrix $A=(\a_{ij})$:
\begin{equation}
\begin{array}{l}
g:=6; \ a:=[1..g]; \ b:=a; \ c:=-a;\\
M:=\operatorname{NullMat}(g,g);\\
\text{ for } i \text{ in } [1..g] \text{ do}\\ 
\text{ for } j \text{ in } [1..g] \text{ do}\\
M[i][j]:= 1/((a[j]-c[i])*(b[j]-c[i]));\text{ od; od};\\
N:=\operatorname{NullMat}(g,g);\\
\text{ for } i \text{ in } [1..g] \text{ do}\\ 
\text{ for } j \text{ in } [1..g] \text{ do}\\
N[i][j]:=(2*c[i]-a[j]-b[j])/((a[j]-c[i])\hat{\phantom{x}}{2}*(b[j]-c[i])\hat{\phantom{x}}{2});\text{ od; od};\\
A:=N/M;\nonumber
\end{array}
\end{equation}

Calculating the number of free variables in the system \eqref{ct-hom-rel},
where we write the coefficients of each equation in a $2g\times g$-matrix, with one block corresponding to
the variables $(\ga_{ij})$ and the other to the variables $(\b_{ij})$:
\begin{equation}
\begin{array}{l}
T := \operatorname{Tuples}([1..g],3);\\
\text{ for } S \text{ in }\operatorname{Tuples}([1..g],3) \text{ do}\\
\text{ if } S[1]>=S[2] \text{ or } S[1]=S[3] \text{ or } S[2]=S[3] \text{ then}\\
\operatorname{RemoveSet}(T, S); \text{ fi; \ od;}\\
\operatorname{equations}:=[];\\
\text{ for } S \text{ in } T \text{ do}\\
m:=\operatorname{NullMat}(2*g,g);\\
m[S[1]+g][S[3]]:=A[S[2]][S[1]]; \ m[S[2]+g][S[3]]:=A[S[1]][S[2]];\\
m[S[2]][S[3]]:=-A[S[1]][S[3]]; \ m[S[2]][S[1]]:=A[S[1]][S[3]];\\
m[S[1]][S[3]]:=-A[S[2]][S[3]]; \ m[S[1]][S[2]]:=A[S[2]][S[3]];\\
\operatorname{Add}(\operatorname{equations},m); \text{ od};\\
V :=  \operatorname{FreeLeftModule}(\operatorname{Rationals}, \operatorname{equations});\\
2*g*(g-1)-\operatorname{Dimension}(V);\nonumber
\end{array}
\end{equation}

Calculating the number of free variables in the system \eqref{acd-hom-rel},
where we write the coefficients of each equation in a $3g\times g$-matrix, with blocks corresponding
to the variables $(\ga_{ii})$, $(a_{ij})$ and $(\varepsilon_{ij})$, respectively:
\begin{equation}
\begin{array}{l}
\operatorname{equations2}:=[];\\
\text{ for } S \text{ in } T \text{ do }\\
m:=\operatorname{NullMat}(3*g,g);\\
m[S[1]+2*g][S[3]]:=A[S[2]][S[1]]; \ m[S[2]+2*g][S[3]]:=A[S[1]][S[2]];\\
m[S[1]+g][S[2]]:=1;\ m[S[3]][S[3]]:=3*A[S[1]][S[3]]*A[S[2]][S[3]];\\
\operatorname{Add}(\operatorname{equations2},m); \text{ od};\\
V :=  \operatorname{FreeLeftModule}(\operatorname{Rationals}, \operatorname{equations2});\\
3*g*(g-1)/2+g-\operatorname{Dimension}(V);\nonumber
\end{array}
\end{equation}

\noindent
{\bf 2. GAP codes for Theorem \ref{dominant-thm}}. 

\medskip

Setting up (say, for genus $5$) and calculating matrices $M$, $N$, $A$ as before,
as well as some auxiliary quantities, namely, the matrices
$ecp=(e'_j(c_i))$, $ecpp=(e''_j(c_i))$,
$fpa=(f'_i(a_j))$, $fpb=(f'_i(b_j))$, $eta=(\wt{e}_{ji}(a_i))$ and $etb=(\wt{e}_{ji}(b_i))$:
\begin{equation}
\begin{array}{l}
g:=5;\  a:=[1..g];\ b:=[(g+1)..(2*g)];\ c:=-a;\\
M:=\operatorname{NullMat}(g,g);\ N:=\operatorname{NullMat}(g,g);\ 
Np:=\operatorname{NullMat}(g,g);\\
ac2:=\operatorname{NullMat}(g,g);\ 
ac3:=\operatorname{NullMat}(g,g);\ 
bc2:=\operatorname{NullMat}(g,g);\ 
bc3:=\operatorname{NullMat}(g,g);\\
\text{ for } i \text{ in } [1..g] \text{ do }\\
\text{ for } j \text{ in } [1..g] \text{ do }\\
M[i][j]:=(a[j]-c[i])\hat{\phantom{x}}(-1)*(b[j]-c[i])\hat{\phantom{x}}(-1);\\
N[i][j]:=(2*c[i]-a[j]-b[j])*(a[j]-c[i])\hat{\phantom{x}}{(-2)}*(b[j]-c[i])\hat{\phantom{x}}(-2);\\
Np[i][j]:=2*((c[i]-a[j])\hat{\phantom{x}}(-3)-(c[i]-b[j])\hat{\phantom{x}}(-3))/(a[j]-b[j]);\\
ac2[i][j]:=(a[j]-c[i])\hat{\phantom{x}}(-2);\
ac3[i][j]:=(a[j]-c[i])\hat{\phantom{x}}(-3);\\
bc2[i][j]:=(b[j]-c[i])\hat{\phantom{x}}(-2);\
bc3[i][j]:=(b[j]-c[i])\hat{\phantom{x}}(-3);\
\text{ od;\ od};\\
A:=N/M;\\
x:=\operatorname{Inverse}(M);\
ecp:=-N*x;\
ecpp:=Np*x;\\
fpa:=2*ac3+A*ac2;\
fpb:=2*bc3+A*bc2;\\
eb:=\operatorname{NullMat}(g,g);\
ea:=\operatorname{NullMat}(g,g);\\
\text{ for } i \text{ in } [1..g] \text{ do }\\
\text{ for } j \text{ in } [1..g] \text{ do }\\
\text{ if } i=j \text{ then } eb[i][i]:=-(b[i]-a[i])\hat{\phantom{x}}(-2);\
ea[i][i]:=eb[i][i];\\
\text{ else } eb[i][j]:=(b[i]-a[j])\hat{\phantom{x}}(-1)*(b[i]-b[j])\hat{\phantom{x}}(-1);\\
ea[i][j]:=(a[i]-a[j])\hat{\phantom{x}}(-1)*(a[i]-b[j])\hat{\phantom{x}}(-1);\
\text{ fi}; \text{ od;\ od};\\
etb:=eb*x;\
eta:=ea*x;\nonumber
\end{array}
\end{equation}

In the main cycle we create the set of vectors of length $5g$ numbered by pairs $(i,j)$, $i\neq j$,
representing images of $\eta_j\in H^0(C,\om_C(-p_i))$ under \eqref{nodal-dual-tangent-eq}.
The coordinates of these vectors are partitioned into $5$ segments of length $g$
(named $tau$, $pa$, $pb$, $pc1$ and $pc2$), corresponding
respectively to 
$\tau_{q_k}(df_i\ot\eta_j)$, and the polar parts of the Laurent expansions of $df_i\ot\eta_j$
at $a_k$, $b_k$ and $c_k$ (the latter are recorded in two segments: 
coefficients of $(dt)^{\ot 2}{t-c_k}$
in positions $[3g+1,\ldots,4g]$ and coefficients of $(dt)^{\ot 2}{(t-c_k)^2}$).
The output is the rank of the map \eqref{nodal-dual-tangent-eq}.

\begin{equation}
\begin{array}{l}
\operatorname{functionals}:=[];
\text{ for } i \text{ in } [1..g] \text{ do }\\
\text{ for } j \text{ in } [1..g] \text{ do }\\
\text{ if } i<>j \text{ then }\\
tau:=0*[1..g];\
pa:=0*[1..g];\
pb:=0*[1..g];\
pc1:=0*[1..g];\
pc2:=0*[1..g];\\
xi:=0*[1..(5*g)];\\
\text{ for } k \text{ in } [1..g] \text{ do }\\
tau[k]:=etb[k][j]*fpa[i][k]+eta[k][j]*fpb[i][k];\\
pa[k]:=x[k][j]*fpa[i][k];\
pb[k]:=x[k][j]*fpb[i][k];\\
\text{ if } k=i \text{ then }
pc1[i]:=ecpp[i][j]+ecp[i][j]*A[i][i];\
pc2[i]:=2*ecp[i][j];\\
\text{ elif } k=j \text{ then }
pc1[j]:=ecp[j][j]*A[i][j];\
pc2[j]:=A[i][j];\\
\text{ else } pc1[k]:=ecp[k][j]*A[i][k];\
\text{ fi};\\
xi[k]:=tau[k];\
xi[k+g]:=pa[k];\
xi[(k+2*g)]:=pb[k];\\
xi[(k+3*g)]:=pc1[k];\
xi[(k+4*g)]:=pc2[k];\
\text{ od };\\
\operatorname{Add}(\operatorname{functionals},xi);\\
\text{ fi};\
\text{ od;\ od};\\
V :=  \operatorname{FreeLeftModule}( \operatorname{Rationals}, \operatorname{functionals});\\
\operatorname{Dimension}(V);\nonumber
\end{array}
\end{equation}

For $g=4$ and $g=5$ we get the rank equal to $12$ and $20$, respectively. As a sanity check,
for $g=6$ we get the rank equal to $27=5g-3$ 
which is the dimension of the moduli space $\MM_{6,6}^{(1)}$,
which agrees with the fact that for $g\ge 6$ the tangent map is generically injective.

\bigskip

{\sc Department of Mathematics, University of Oregon, Eugene, OR 97403}

{\it Email addresses}: rfisette@@uoregon.edu, apolish@@uoregon.edu


\begin{thebibliography}{99}
\bibitem{Bardzell} M.~J.~Bardzell, {\it The alternating syzygy behavior of monomial algebras},
J.~Algebra 188 (1997), 69--89.
\bibitem{BK} A.~Bondal, M.~Kapranov,
{\it Framed triangulated categories}, Math. USSR-Sb.~70 (1991), 93--107.
\bibitem{Cohen} J.~M.~Cohen, {\it The decomposition of stable homotopy}, Annals Math.
87 (1968), 305--320.
\bibitem{Efimov} A.~I.~Efimov,
{\it Homological mirror symmetry for curves of higher genus}, Advances Math. 230 (2012), 493--530.
\bibitem{Fisette-thesis} R.~Fisette, {The A-infinity algebra of a curve and the $j$-invariant}, Ph.D. thesis, 
University of Oregon, 2012.
\bibitem{GM}  S.~Gelfand, Yu.~Manin, {\it Methods of homological algebra}. Springer-Verlag, 1996.
\bibitem{GL} M.~Green, R.~Lazarsfeld, {\it Some results on the syzygies of finite sets and algebraic curves},
Compositio Math. 67 (1988), 301--314.
\bibitem{GS} V.~K.~A.~M.~Gugenheim, J.~D.~Stasheff, {\it On
perturbations and $A\sb \infty$-structures}. Bull. Soc. Math. Belg.
   Sir. A 38 (1986), 237--246. 
\bibitem{GLS} V.~K.~A.~M.~Gugenheim, L.~A.~Lambe, J.~D.~Stasheff, 
{\it Perturbation theory in differential homological
   algebra. II.} Illinois J. Math. 35 (1991), no. 3, 357--373.
\bibitem{Ha-DT} R.~Hartshorne, {\it Deformation Theory}, Springer-Verlag, New York, 2010.
\bibitem{HVdB} L.~Hille, M.~Van den Bergh, {\it Fourier-Mukai transforms}, in {\it Handbook of tilting theory}, 
147--177, Cambridge Univ. Press, Cambridge, 2007.
\bibitem{Kadeishvili}
T.~V.~Kadeishvili, {\it The category of differential
coalgebras and the category of $A_{\infty}$-algebras} (in Russian).
Trudy Tbiliss. Mat. Instituta 77 (1985), 50--70.
\bibitem{Kap} M.~Kapranov, {\it Derived categories of coherent sheaves on homogeneous spaces},
dissertation (in Russian), Steklov Math. Institute, 1988.
\bibitem{KKOY} A.~Kapustin, L.~Katzarkov, D.~Orlov, M.~Yotov, {\it Homological mirror symmetry for manifolds of general type}, Cent. Eur. J. Math. 7 (2009), 571--605.
\bibitem{Keller-intro} B.~Keller, {\it Introduction to A-infinity algebras and modules}, Homology Homotopy Appl. 3 (2001), 1--35.
\bibitem{Keller} B.~Keller, {\it $A$-infinity algebras, modules and functor categories}, in {\it Trends in representation theory of algebras and related topics}, 67--93, AMS, Providence, RI, 2006.
\bibitem{KS} M.~Kontsevich, Y.~Soibelman, {\it Homological mirror symmetry
and torus fibration}, in {\it Symplectic geometry and mirror 
symmetry (Seoul, 2000)}, 203--263, World Sci. Publishing, River Edge, NJ, 2001. 
\bibitem{Kraines} D.~A.~Kraines, {\it Higher Massey products}, Trans.~AMS 124 (1966), 431--439.
\bibitem{LP} Y.~Lekili, T.~Perutz, {\it Fukaya categories of the torus and Dehn surgery}, Proc. Natl. Acad. Sci.
108 (2011), 8106--8113.
\bibitem{LPWZ} D.-M. Lu, J. H. Palmieri, Q.-S. Wu and J. J. Zhang, {\it A-infinity structure on Ext-algebras}, J. Pure Appl. Algebra 213 (2009) 2017--2037,
\bibitem{May} J.~P.~May, {\it Matric Massey Products}, J. Algebra 12 (1969), 533-568.
\bibitem{Merkulov} S.~Merkulov, {\it Strong homotopy algebras of a K\"ahler
manifold}, Internat. Math. Res. Notices 1999, no.3, 153--164.
\bibitem{Orlov} D.~Orlov, {\it Remarks on generators and dimensions of triangulated categories}, 
Moscow Math. J. 9 (2009), 153--159.
\bibitem{P-ext} A.~Polishchuk, {\it Extensions of homogeneous coordinate rings to
          $A_{\infty}$-algebras}, Homology Homotopy Appl.          5 (2003), 407--421. 
\bibitem{P-Mas} A.~Polishchuk,
{\it Triple Massey products on curves, Fay's trisecant identity
          and tangents to the canonical embedding}, Moscow Math. J. 3 (2003), 105--121. 
\bibitem{P-ell} A.~Polishchuk, {\it $A_{\infty}$-algebra of an elliptic curve and Eisenstein series}, 
Comm. Math. Phys. 301 (2011), 709--722.
\bibitem{Seidel} P.~Seidel, {\it Homological mirror symmetry for the genus two curve}, J. Algebraic Geom. 20 (2011), no. 4, 727--769.
\bibitem{Shipley} B.~Shipley, {\it An algebraic model for rational $S^1$-equivariant stable homotopy theory},
Q. J. Math. 53 (2002), 87--110.
\bibitem{Toen} B.~To\"en, {\it Finitude homotopique des dg-algebres propres et lisses},
Proc. Lond. Math. Soc. (3) 98 (2009), 217--240.
\bibitem{Wahl} J.~Wahl, {\it Gaussian maps on algebraic curves}, J. Differential Geom. 32 (1990), 77--98.
\end{thebibliography}
\end{document}